\documentclass[a4paper,twoside,10pt]{article}

\usepackage{amssymb}
\usepackage{amsmath}
\usepackage{epsfig}

\newtheorem{The}{Theorem}
\newtheorem{Lem}{Lemma}

\newtheorem{Pro}[Lem]{Proposition}

\newtheorem{Rem}[Lem]{Remark}
\newtheorem{Defi}[Lem]{Definition}

\newtheorem{Not}[Lem]{Notation}

\def\a{\alpha}

\def\C{\mathbf C}
\def\D{\Delta}
\def\d{\delta}
\def\e{\varepsilon}
\def\f{\phi}
\def\g{\gamma}

\def\l{\lambda}

\def\O{\mathbf O}
\def\p{\pi}

\def\Q{\mathbf Q}
\def\R{\mathbf R}
\def\r{\rho}

\def\s{\sigma}
\def\t{\tau}

\def\w{\omega}
\def\Z{\mathbf Z}
\def\z{\zeta}

\def\elem(#1,#2){  {{#1}\over \overline {\ #2\ }} }  

\title{Decomposition in bunches \\ of the critical locus of a quasi-ordinary map}
\author{E.R. Garc\'\i a Barroso, P.D. Gonz\'alez P\'erez}
\date{\,}
\begin{document}

\maketitle

\begin{quote} {\bf Abstract.} 
\footnote{ Math. classification numbers: 14M25, 32S25.

} 
A polar hypersurface $P$ of a complex analytic hypersurface germ $f= 0$ can be 
investigated 
by analyzing the invariance of certain Newton polyhedra 
associated to the  image of $P$, with respect to suitable coordinates, 
by certain
morphisms appropriately associated to $f$.
We develop this general principle of Teissier (see  \cite{TP} and \cite{TJ})
when $f= 0$ is a quasi-ordinary hypersurface germ and $P$ is the polar hypersurface associated to 
any quasi-ordinary projection of $f=0$. 
We build a  decomposition of $P$ in
bunches of branches which characterizes 
the embedded topological type of the irreducible components of $f= 0$.
This decomposition is characterized also by 
some properties of the strict transform of $P$ by the toric embedded  resolution of $f = 0$ given by the second author in \cite{GP3}.
In the plane curve case this result provides a simple algebraic proof of a theorem of L\^e, Michel and Weber in \cite{LMW}.
\end{quote}

\section*{Introduction}

The {\it polar varieties} or at least their rational equivalence classes play 
an important role in projective geometry in particular 
in the study of characteristic classes and numerical invariants of projective algebraic varieties, 
and also in the study of projective duality 
(Pl\" ucker formulas). In the 1970's local polar varieties began to be used systematically in the 
study of singularities. Local polar varieties can be used to produce 
invariants of equisingularity (``topological'' invariants of complex analytic singularities)
and also to explain why the same invariants appear in apparently unrelated questions. 
We study here a particular instance of construction and study of such equisingularity invariants. 

The {\it Jacobian polygon}, a plane polygon associated by Teissier 
to a germ of complex analytical hypersurface  
defining an isolated singularity at the origin, 
is an invariant of equisingularity for {\it $c$-equinsingularity},
which is equivalent to Whitney condition for a family of isolated hypersurface singularities, 
implies topological triviality and is equivalent to it for plane curves. 
The inclinations of the compact edges of this polygon are rational numbers called 
the {\it polar invariants} of the germ. 
If   $\f(X_1, \dots, X_{d+1}) =0$ is the equation of the hypersurface germ 
with respect to some suitable coordinates 
the Jacobian polygon coincides with the Newton 
polyhedron of image of the critical locus, or {\it  polar variety}, 
of the morphism $(\C^{d+1}, 0) \rightarrow (\C^2, 0)$
defined by $T = \f(X_1, \dots, X_{d+1})$ and  $U = X_1$,
with respect to the coordinates $(U, T)$ (see \cite{TP}).

In the case of a germ of plane irreducible curve, 
Merle shows in \cite{Merle2} that the polar invariants 
determine also the equisingularity class of the curve 
(or equivalently its embedded topological type).
Merle's results has been generalized to the case of {\it reduced} plane curve germs 
by Kuo, Lu, Eggers, Garc\' \i a Barroso and Wall among others (see \cite{KL}, \cite{Eggers}, 
\cite{GB} and \cite{Wall}). 
They give a decomposition theorem of a generic polar curve 
in bunches which depends only on the equisingularity class of the curve. 
To this  decomposition is associated a matrix of
{\it partial polar invariants}  
which determines the equisingularity 
class of the curve (see \cite{GB})).
L\^e, Michel and Weber 
have proven using topological methods that 
the strict transform of a generic polar curve 
by the minimal embedded resolution of the curve
meets any connected  component of a {\it permitted} subset of 
the exceptional divisor (see \cite{LMW}). 
Another decomposition in bunches of the generic polar curve, 
can be defined from this result in a 
geometrical way, the bunches correspond 
bijectively to the connected components of the permitted subset.
Garc\'\i a Barroso has compared these two decompositions 
and shown that they coincide in  \cite{GB}.

In this paper we study local polar hypersurfaces of a class of 
complex analytic hypersurface singularity, called {\it quasi-ordinary}.
This class of singularities, of which the simplest example are the singularities of plane curves, 
appears naturally in 
Jung's approach to analyze surface  singularities and their parametrizations.
A germ  of  complex analytic variety
is {\it quasi-ordinary} if there exists
a finite projection, called quasi-ordinary,
to the complex affine space  
with discriminant locus contained in  a normal crossing divisor.
By Jung-Abhyankar's theorem 
any quasi-ordinary projection is provided with 
a
parametrization with fractional power series (\cite{Jung} and \cite{Abhyankar}).
In the hypersurface case these parametrisations determine
a finite set of monomials, called {\it characteristic} or {\it distinguished},
which determine quite a lot of the geometry and topology of the singularity.
For instance,  
these monomials
constitute a complete invariant of its {\it embedded topological type}  
in the analytically irreducible case 
and conjecturally in the reduced case (see Gau and Lipman works \cite{Gau}, \cite{Lipman1} and 
\cite{Lipman2}), in particular they determine the zeta function of the geometric monodromy 
as shown in the works of N\'emethi, McEwan and Gonz\'alez P\'erez (see \cite{zeta1} and \cite{zeta2}).
The characteristic monomials 
determine also embedded resolutions of the corresponding 
quasi-ordinary hypersurface singularity, which have been obtained 
in two different ways by Villamayor \cite{Villamayor}
and Gonz\'alez P\'erez (see \cite{GP3} and \cite{GP4}).

We give a decomposition theorem of 
the polar hypersurfaces  $(P, 0)$ of a quasi-ordinary 
hypersurface  $(S, 0)$ corresponding to a 
quasi-ordinary projection.
If $(S, 0)$ is embedded in   $(\C^{d+1} , 0)$
any quasi-ordinary projection can be expressed in
suitable coordinates by  $(X_1, \dots, X_d, Y) \mapsto (X_1, \dots, X_d)$.
Then $(S, 0)$ has an equation defined by a Weierstrass polynomial $ f \in \C \{ X_1, \dots, X_d \} [Y]$
and the associated polar hypersurface $(P, 0)$ is defined as 
the critical space, $f_Y = 0$,      
of the
{\it quasi-ordinary morphism}: 
\[
\left\{
\begin{array}{l}
\xi_f \;\; : \;\; (\C^{d+1}, 0)  \;\; \longrightarrow \;\;  (\C^{d+1}, 0)
\\
U_1 = X_1, \; \; \dots , \;\; U_d = X_d,  \; \; T = f(X_1, \dots,X_d, Y)
\end{array}
\right.
\]
The decomposition is defined in terms of a matrix, 
generalizing the matrix of  partial polar invariants of \cite{GB}, 
which determines and is determined by the 
partially ordered set of {\it characteristic monomials} associated to the
fixed quasi-ordinary projection. 
In particular, it defines a complete invariant of the embedded topological type
of each irreducible component of  $(S, 0)$ by using Gau's and Lipman's results.
Our decomposition theorem is partially motivated by 
a result of Popescu-Pampu providing a decomposition of the polar hypersurface 
$(P, 0)$
of the quasi-ordinary hypersurface $(S, 0)$, 
obtained
with the additional hypothesis that $(S, 0)$ and $(P, 0)$ are 
simultaneously quasi-ordinary with respect to the given
quasi-ordinary projection (see  \cite{PP2}, Chapter 3 or \cite{PP1}). 
Our decomposition 
extends also to the {\it Laurent} quasi-ordinary case
 studied in \cite{PP2} by analogy to the case of {\it meromorphic} plane curves
of Abhyankar and Assi (see \cite{A-A}).

An important tool which we introduce  to prove these results is 
 an {\it irreducibility criterion} for 
power series with {\it polygonal Newton polyhedra} 
(the maximal dimension of its compact faces is equal to one).
Our  criterion, which holds for power series
with coefficients over any algebraically closed field of arbitrary 
characteristic, 
states that if an irreducible series has a polygonal Newton 
polyhedron then it has only one compact edge. 
This result generalizes a fundamental property of plane curve germs.
Our proof is obtained by using Newton polyhedra in the 
framework of toric geometry.

The decomposition theorem of the polar hypersurface $(P, 0)$ 
has a proof inspired by Teissier's works \cite{TP} and
\cite{TJ}. In the irreducible case we analyze the 
{\it discriminant} ${\cal D}$ of the quasi-ordinary map $\xi_f$, i.e., 
the   image of the critical space.
We compute  the Newton polyhedron  ${\cal N }_{\cal D}$ of 
the discriminant  ${\cal D}$ in the coordinates $U_1, \dots, U_d, T$ above,
and we show that it is a {\it polygonal polyhedron}.
This computation applies the above mentionned  results of Popescu-Pampu,  after suitable toric base changes
already used in \cite{GP1}.
We use then the irreducibility criterion 
to show the existence of a 
decomposition of $(P, 0)$ in bunches 
which 
correspond bijectively to the compact edges of the  polyhedron ${\cal N }_{\cal D}$.

We also give  a  geometrical characterization of the decomposition theorem 
by analyzing  the strict transform of $(P, 0)$ 
by a modification $p: {\cal Z} \rightarrow \C^{d+1}$
which is built canonically from the given
quasi-ordinary projection by using the characteristic monomials (see \cite{GP4}). 
Geometrically the bunches of the decomposition of $(P,0)$
correspond to the union of branches  of the polar hypersurface $(P, 0)$  
whose strict transforms by $p$,
meet the  exceptional fiber  $p^{-1} (0) $  at the same irreducible component.
A posteriori this analysis can be extended to the
the toric {\it embedded resolutions} of $(S,0)$ built in \cite{GP3} or \cite{GP4}, since 
they are  factored by $p$. 
In the plane curve case we apply this result to obtain a simple algebraic proof of the theorem of 
L\^e, Michel and Weber in \cite{LMW} which shows the underlying 
toric structure of the decomposition of the polar curve.

Our results provide answers to some of the questions 
raised independently by McEwan and N\'emethi in \cite{RM} section III, among 
some open problems concerning quasi-ordinary singularities.
We hope that the results of this paper could apply to the study of polar varieties of 
hypersurface singularities by using a suitable form of Jung's approach.
It is reasonable to expect that this work may have some applications to the metric study of 
the Milnor fibres of hypersurfaces, at least in the quasi-ordinary case, 
as suggested by Teissier's and Garc\'\i a Barroso's results 
in the case of  plane curve singularities (see  \cite{Barroso}); 
see also Risler's work \cite{Risler} for the real plane curve case. 

The proofs are written in the analytic case. The results and proofs 
of this paper hold also  in the algebroid case
(over an algebraically closed field of  characteristic zero).

{\bf Acknowledgments.} 
We are grateful to Bernard Teissier for his suggestions and comments.
This research has been partly financed by ``Acci\'on integrada hispano-francesa HF 2000-0119'' and  by ``Programme d'actions int\'egr\'ees franco-espagnol
02685ND''.
The second author is supported by a Marie Curie Fellowship of the European Community
program ``Improving Human Research Potential and the Socio-economic Knowledge Base'' 
under contract number HPMF-CT-2000-00877.

\section{Quasi-ordinary polynomials, their characteristic monomials and the Eggers-Wall tree}
A germ of complex analytic hypersurface $(S, 0) \subset  (\C^{d+1}, 0) $   
is {\it quasi-ordinary} if there exists a finite projection 
$(S, 0) \rightarrow (\C^d,0)$
which is a local isomorphism outside a normal crossing divisor. 
The embedding $(S, 0) \subset  (\C^{d+1}, 0) $   can be defined by an equation 
$f= 0$ where $f \in \C \{ X \} [Y]$ is 
a {\it quasi-ordinary polynomial}: a  Weierstrass polynomial 
with discriminant $\D_Y f$ of the form
$\D_Y f = X^\d \epsilon$ for a unit $\epsilon$ in the ring 
$ \C \{ X \}$ of convergent (or formal) power series in the variables 
$X= (X_1, \dots, X_d)$ and $\d \in \Z^d_{\geq 0}$.

The Jung-Abhyankar theorem 
guarantees that the roots of the quasi-ordinary polynomial $f$ 
are fractional power series in the ring $\C \{ X^{1/k} \}$ for some suitable integer $k$,
for instance $k =\deg f$ when $f$ is irreducible
(see \cite{Abhyankar}).
If the series 
$\{ \z ^{(l)} \}_{l =1}^{\deg f}  \subset \C \{ X^{1/k} \}$ are 
the roots of $f$, its discriminant is equal to: 
\begin{equation} \label{disc}
\D_Y f = \prod_{i\ne j} (\z^{(i)} - \z^{(j)})
\end{equation}
hence 
each factor $\z^{(t)} - \z^{(r)}$ is  of the form 
$X^{\l_{t,r}}  \epsilon _{t,r} $ where $\epsilon_{t,r} $ is a unit in $\C \{ X^{1/k} \} $.
The monomials $X^{\l_{t,r}}$ 
(resp. the exponents ${\l_{t,r}}$) 
are called {\it characteristic}.

In the reducible reduced case,
if $f=  f_1 \dots f_s $ is the factorization in monic irreducible polynomials
each factor $f_i$ is a quasi-ordinary polynomial 
since $\D_Y f_i$ divides $\D_Y f$ by formula (\ref{disc}). 

We define the partial order 
of $\R^d \cup \{ + \infty \}$: 
\begin{equation} \label{order}
u \leq u' \Leftrightarrow u' \in u + \R_{\geq 0}^d.
\end{equation}
We write $u < u'$ if $u \leq u'$ and $u \ne u'$. If $\a \in \R^d$ we set $ \a < + \infty$. 
Notice that $u \leq u' $ means that the inequality holds 
coordinate-wise with respect to the canonical 
basis.
The characteristic exponents have the following property with respect to the order (\ref{order}), 
see \cite{Lipman2} and \cite{Licei}.

\begin{Lem} \label{rorder} 
Let $f_i$ be an irreducible factor of the reduced quasi-ordinary polynomial $f$.
The set 
\begin{equation} \label{theta}
V_f (f_i) := 
\left\{ \l_{r,t} / \z^{(r)} \ne \z^{(t)}, \; f (\z^{(t)}) = 0 \mbox{ and } f_i (\z^{(r)}) = 0 \right\}
\end{equation}
is totally ordered by  $\leq$.
\hfill $\ {\diamond}$
\end{Lem}

If $f_i$ and $f_j$ are two  irreducible factors of the quasi-ordinary polynomial $f$
we define the {\it order of coincidence}  $k(f_i, f_j)$  of their roots by: 
\[
k(f_i, f_j) = \max \{ \l_{r,t}  /   f_i (\z^{(r)}) = 0,    f_j (\z^{(t)}) = 0  \}
\]
The order of  coincidence of $f_i $ with itself is $k(f_i, f_i):=+ \infty$.
We have the following ``valuative'' property of the orders of coincidence (see Lemma 3.10 of \cite{GP3}). 
\begin{equation} \label{val}
\min\{ k(f_i, f_j) , k(f_j, f_r) \} \geq k(f_i, f_r) \mbox{ with equality if }  k(f_i, f_j) \ne k(f_j, f_r)
\end{equation}

The totally ordered set $V_f (f_i)$ defined by (\ref{theta}) is equal to the union of the non necessarily disjoint sets
whose elements are 
the characteristic exponents
 $\l_1^{(i)} <  \cdots <  \l_{g(i)} ^{(i)} $ of $f_i$, if they 
exist\footnote{The case $f_i$ with no characteristic monomials happens only when $\deg f_i =1$}, and  
the orders of coincidence $ k(f_i, f_j)$ for $j= 1, \dots , s$ and $j\ne i$.
We associate to
the characteristic exponents of the irreducible factor $f_i$, for $i=1, \dots, s$,
the following sequences of 
{\it characteristic lattices and integers}:
the lattices are $M_0^{(i)} :=\Z^{d} $  and  $M_j^{(i)} := M_{j-1}^{(i)} + \Z \l_j^{(i)}$
for $j=1, \dots, g(i)$ with the convention $\l_{g(i) +1}^{(i)} = + \infty$;
the integers are $n_0^{(i)}:=1$ and $n_j^{(i)}$ is the index of the subgroup $ M_{j-1}^{(i)}$
in $M_j^{(i)}$,  for $j=1, \dots, g(i)$.
We denote the integer $n_j^{(i)}  \cdots n_{g(i)}^{(i)}$ by
$e_{j-1}^{(i)}$ for $j= 1, \dots , g(i)$. 
We have that $\deg f_i =  e_0^{(i)}=  n_1^{(i)}  \cdots n_{g(i)}^{(i)}$
(see \cite{Lipman2} and \cite{Tesis})
and used
When $d=1$  we have the equality  $M_{j}^{(i)} = (e_j^{(i)}) \Z$, 
and the integer  $n_j$ coincides with 
the first component of
the classical characteristic pairs of the plane branch defined by $f=0$.

The information provided by the characteristic monomials 
is structured in a tree which encodes the embedded topological type 
of the irreducible components of $f= 0$ (characterized by the work of 
Gau and Lipman in terms of the characteristic exponents, see \cite{Gau} and \cite{Lipman2}).
This tree is introduced by
Popescu-Pampu (see \cite{PP2} and 
\cite{PP1}) following a construction of Wall \cite{Wall} and Eggers \cite{Eggers}
used to study the polar curves of a plane curve germ (see also \cite{GB}).

The {\it elementary branch} $ \theta_f (f_i)$ associated to $f_i$ 
is the abstract simplicial complex of dimension one 
with vertices running through 
the elements of the totally ordered subset 
$ V_{f} (f_i) \cup \{ 0, +  \infty \} $ of $\Q^d \cup  \{       \infty \}$,
and edges running through the segments joining 
consecutive vertices. 
The underlying topological space is homeomorphic to the segment $[0, + \infty ]$.
We denote the vertex of  $\theta_f (f_i)$ corresponding to   $\l \in 
V_{f} (f_i) \cup \{ 0, +  \infty \} $
by $P^{(i)}_\l $.
The simplicial complex $ \theta_f (f)$  obtained from the disjoint union 
$\bigsqcup_{i=1}^{s} \theta_f (f_i)$
by identifying  in $ \theta_f (f_i)$ and $ \theta_f (f_j)$ the sub-simplicial complexes 
corresponding to $ \overline{ P^{(i)}_{0} P^{(i)}_{k(f_i, f_j)}} $ and 
$ \overline{ P^{(j)}_{0} P^{(j)}_{k(f_i, f_j)} }$ for $1 \leq i <j \leq s$ is a tree.

We give to a vertex 
$P_\l^{(i)} $ of $ \theta_f (f)$ the valuation  $v (P_\l^{(i)}) =  \l$.
The restriction of the valuation $v$ to the set of non extremal vertices of 
$ \theta_f (f)$ is a $0$-chain 
with coefficients in $\frac{1}{k} \Z^d$.
The set of vertices of $\theta_f (f) $ is partially ordered by $P \leq  P'$ if 
$P, P'$ are vertices of the same elementary branch of the tree and 
$v (P)  \leq v (P')$.
The valuation $v$  defines an orientation on the tree $\theta_f (f)$.
The boundary operator $\partial$ 
is the linear map of integral 1-chains defined on the segments by 
$\partial ( \overline{PP'} ) = P'- P$ if $v (P) < v (P')$.

If $v(P) \ne +\infty $ the value $v(P)$ is the $d$-uple of coordinates of 
an element of the lattice $\frac{1}{k} \Z^d$, which we denote by $\tilde{v} (P) $,
with respect to the canonical basis of   $\Z^d$. 
If $v(P) = +\infty $ we set $\tilde{v} (P) =  +\infty$.
This defines a {\it lattice valuation} $\tilde{v} $ 
the vertices of $\theta_f(f) $ which is 
preserved by certain modifications (see section \ref{tbc}).
We recover the valuation $v$ from the lattice valuation $\tilde{v}$ and 
the {\it reference lattice cone:}  $( \R^d_{\geq 0}, \Z^d)$
associated to the exponents of the monomials of $\C \{ X \}$.

For $i = 1, \dots, s$ we
define an integral $1$-chain
whose segments are obtained by   
by subdividing the segments of the chain 
\begin{equation} \label{segment}
   \overline{ P_0^{(i)} P^{(i)}_{\l_1^{(i)}}} + 
 n^{(i)}_1 \overline{ P^{(i)}_{\l_1^{(i)}}   P^{(i)}_{\l_2^{(i)}}}   + \cdots + 
n^{(i)}_1 \cdots n^{(i)}_{g(i)} \overline{ P^{(i)}_{ \l_{g(i)}^{(i)} } P^{(i)}_{ + \infty }}
\end{equation}
with the points  corresponding to the orders of coincidence of $f_i$, 
the coefficient of an oriented  segment in the subdivision 
is the same as the coefficient of the oriented segment of (\ref{segment}) containing it.
It follows  that these
$1$-chains  
paste on  $ \theta_f (f)$ defining a $1$-chain which we denote by $\g_f$.
The vertex $P_\l ^{(i)}$, 
if  $\l \ne    0, +\infty $ does  not correspond to  a characteristic exponent of $f_i$
if and only if the vertex $P_\l ^{(i)}$ appears in two segments of $\theta_f (f_i)$ 
with the same coefficient.

\begin{Defi}
The Eggers-Wall tree  of the quasi-ordinary polynomial $f$ 
is the simplicial complex $ \theta_f (f)$
with the chains $ \g_f$ and $v$. 
We denote by  $ \tilde{\theta}_f (f)$ 
the simplicial complex $ \theta_f (f)$
with the chains $ \g_f$ and $\tilde{v}$. 
\end{Defi}

As shown by the work of Wall \cite{Wall} and Popescu-Pampu \cite{PP1},
the Eggers-Wall tree is useful to represent the information provided by the orders of 
coincidence of the roots of $f$ with the roots of $h$ in the set ${\cal RC}_f $ of polynomials
{\it radically comparable} with the polynomial $f$:
\[
{\cal RC}_f := \left\{ h \in \C \{ X \} [Y] \, / h \mbox{ monic and the product  } \, fh  \mbox{ is quasi-ordinary  }
\right\}.
\]
Any $h \in {\cal RC}_f $ is a quasi-ordinary polynomial and the difference of 
its roots with those of $f$ has a dominant monomial 
(viewed in $\C \{ X^{1/k}  \}$ for some suitable $k$).
If $h \in {\cal RC}_f $ we consider the sub-tree $\theta_{fh} (f) = \bigcup_{i=1}^s 
\theta_{fh} (f_i)$
of   $\theta_{fh} (fh)$ as a subdivision of $\theta_f (f)$
induced by $h $. 
If $h$ is irreducible the point $P^{(i)}_{k(f_i, h)}$  is the point of bifurcation of the
elementary  branches 
$\theta_{fh} (f_i)$ and  $\theta_{fh} (h)$ in $\theta_{fh} (fh)$; we denote by $P^{h}_{k(h, f)}$
the point of bifurcation of the elementary branch $\theta_{fh} (h)$ from the 
tree $\theta_{fh} (f)$.
If $h \in {\cal RC}_f$
and if  $h= h_1 \dots h_t $ is the factorization  of $h$ as a product of  monic irreducible polynomials,  
the {\it contact chain} $[h]^{(f)} $ is the integral 
$0$-chain on $\theta_{fh}(f)$ 
defined by: 
$$[h]^{(f)}  = \sum_{j=1}^t \deg h_j P^{h_j} _{k(h_j, f)}. $$
The contact chain $[h]^{(f)}$ is associated with the decomposition 
$h = b_1 \dots b_{s(f,h)}$ in the ring $\C \{ X \} [Y]$, 
where the factors $b_j$  are the products of 
those irreducible factors of $h$ having the same order of 
coincidence  with each irreducible factor of $f$.
\section{Decomposition in bunches and Newton polyhedra of images}

We show that a quasi-ordinary polynomial $f \in \C \{ X \} [Y]$ 
defines in a natural way a decomposition in bunches  
for certain class of polynomials 
which contains the derivative $f_Y$ of $f$.
We state the characterization of the decomposition 
in bunches of the polar hypersurface $f_Y= 0$ and of the 
Newton polyhedron of the equation defining the image of $f_Y= 0$
by the quasi-ordinary morphism $\xi_f$, in terms of the tree $\theta_f(f)$.  

If $f \in \C \{ X \} [Y]$ is a quasi-ordinary polynomial 
the derivative $f_Y := \frac{1}{n} \frac{\partial f}{\partial Y}$
belongs to the set ${\cal C}_f$ of 
 polynomials  {\it comparable}\footnote{Our notion of {\it comparable} polynomials generalizes that of  {\it radically comparable}, 
(which corresponds to the notion of {\em comparable} polynomials of Popescu-Pampu in  \cite{PP2} 
or \cite{PP1}).}  
with $f$:
\[
{\cal C}_f := \left\{ h \in \C \{ X \} [Y] \, / h \mbox{ monic, } \, \mbox{ Res}_Y (f, h) = X^{\r(f, h)} \epsilon_{f, h} \mbox{ with } \epsilon_{f, h} (0) \ne 0  \mbox{ and } \r(f, h) \in \Z^d \right\},
\]
where $\mbox{Res}_Y (f, h)$ denotes the resultant of the polynomials $f$ and $h$.

We have an inclusion ${\cal RC}_f  \subset {\cal C}_f$:
it is sufficient to notice that  $\mbox{Res}_Y (f, h)$ divides $\D_Y (fh)$, 
a statement which follows from the 
classical properties of resultants and discriminants (see \cite{GKZ}): 
\begin{equation} \label{Res}
\mbox{Res}_Y (f, h_1 \cdots h_t ) = \prod_{i=1}^t \mbox{Res}_Y (f, h_i), \quad
\D_Y (fh) = \D_Y (f) \D_Y (h) (\mbox{Res}_Y (f, h))^2
\end{equation}

\begin{Rem}\label{ple}
The inclusion ${\cal RC}_f  \subset {\cal C}_f$ is strict in general.
In particular the derivative $f_Y$ of a quasi-ordinary polynomial $f \in \C \{ X_1, \dots, X_d \}[Y]$ is not 
radically comparable to $f$ in general.
\end{Rem}
- For instance, if $f = Y$ we have that ${\cal C}_f$ is the set 
of monic polynomials $h \in \C \{ X \}[Y]$ such that $h(0) $ is of the form 
$h(0) =  X^{\r(f, h)} \epsilon_{f, h}$ with $\epsilon_{f, h} (0) \ne 0$.
On the other hand  ${\cal RC}_f$ is the set of monic polynomials $h \in \C \{ X \}[Y]$ 
such that the product $Y h$ is a quasi-ordinary polynomial. We have that 
$h = Y^2 + (X_1 + X_2) Y + X_1 X_2^2 \in {\cal C}_f \backslash {\cal RC}_f$.

\noindent
- This exemple is already given by Popescu-Pampu. The polynomial 
$f = Y^{3} + X_1 X_2 Y^2 + X_1^3 X_2 Y + X_1 X_2$ 
is quasi-ordinary thus $f_Y \in {\cal C}_f$  but 
$f_Y \notin {\cal RC}_f$ (see  \cite{PP2}, page 127).

\medskip

If  $f= f_1 \cdots f_s$ is the factorization of $f$ in 
monic irreducible polynomials 
we deduce from  (\ref{Res}) and the definitions that ${\cal C}_f = \bigcap_{i=1}^s {\cal C}_{f_i}$.
We define an equivalence relation in the set 
$  {\cal C}_f $:
\begin{equation} \label{eq}
h, h' \in {\cal C}_f, \,   h \sim h' \Leftrightarrow \frac{\r(f_i, h)}{\deg h} =  
\frac{\r(f_i, h')}{\deg h'}
\mbox{ for } i=1, \dots, s. 
\end{equation}
By (\ref{Res})  if $h \in {\cal C}_f$ the irreducible factors of $h$ are also in ${\cal C}_f$.
We denote by  $s(f,h)$ the number of  
classes 
of the restriction of the equivalence relation (\ref{eq})
to the set of  irreducible factors of $h$. 

\begin{Defi} \label{type}
The $f$-bunch  decomposition of a polynomial $h \in \C \{ X_1, \dots, X_d \} [Y]$
comparable to the quasi-ordinary polynomial $f \in \C \{ X_1, \dots, X_d \} [Y]  $ is 
$
h = b_1 \cdots b_{s(f,h)}
$, 
where the $b_i$, for $i=1, \dots, s(f,h)$, are the products of the irreducible factors 
of $h$ which are in the same class. 
The {\it type of the $f$-bunch decomposition of $h$}  is the collection of 
vectors: 
\begin{equation}
\left\{  \left( \frac{\r(f_1, h_j)}{\deg{h_j}}, \dots,
\frac{\r(f_s, h_j)}{\deg{h_j}};   \deg {b_j} \right) \right\}_{j=1}^{s(f,h)},
\end{equation}
where $h_j$ is any irreducible factor of $b_j$ for $j =1, , \dots, s(f,h)$.
\end{Defi}
If $f$ is clear from the context, we write  
bunch decomposition instead of $f$-bunch decomposition, in particular 
we will do this for the polynomial  $f_Y$.

If $h \in {\cal RC}_f$ the type of the $f$-bunch decomposition of $h$ is studied by 
using the orders of coincidence of $h \in {\cal RC}_f$ with the irreducible factors of $f$.
The following proposition extends the classical relation between 
the intersection multiplicity and the order of coincidence 
in the plane branch case (see Proposition 3.7.15 of  \cite{PP2}).
We introduce some notations: 
if $\l \in V_{hf} (f_i)$ we denote by 
$c_{\l} ^{(i)}$  the integer $c_{\l} ^{(i)} = \max ( \{j / \l^{(i)}_j < \l \} \cup \{ 0 \} )$, 
if $ \l = k(h, f_i) $ we denote $c_{\l} ^{(i)}$ also by $c_{(h, f_i)}$.

\begin{Pro} {\rm  \label{coincidence}}
Let $h \in {\cal RC}_f$ irreducible. If $\t$ is any root of $h$ in 
$\C \{ X^{1/k} \}$ (for some integer $k > 0$) we have that $f_i (\t) $ is of the form:
\begin{equation}  \label{froots}
f_i( \t) = X^{ \frac{ \r (f_i, h)  } {\deg h} } \epsilon_{i,\t} \mbox{ where } \epsilon_{i,\t} 
\mbox{ is a unit in  } \C \{ X^{1/k} \} \mbox{ and }
\end{equation}
\begin{equation} \label{19}
\frac{ \r (f_i, h) } {\deg h} = e_{c_{(h, f_i)}}^{(i)}  k(h, f_i) +  \sum_{k=1}^{c_{(h, f_i)}} (e_{k-1}^{(i)} - e_k^{(i)} ) \l_k^{(i)}
\end{equation}
\hfill $\ {\diamond}$
\end{Pro} 


\begin{Defi} 
We associate to the factor $f_i$ of $f$ the valuation $\nu_i$ of the vertices
of $\theta_{fh} (f)$:
\begin{equation} \label{nu}
\nu_i (P_\l ^{(j)} ):= \left\{ \begin{array}{lll}
e_{c_{\l} ^{(i)}}^{(i)}  \l +  
\sum_{k=1}^{c_{\l} ^{(i)}} (e_{k-1}^{(i)} - e_k^{(i)} ) \l_k^{(i)} & \mbox{ if } & P_\l ^{(j)} \in \theta_{fh} (f_i), \l \ne 0 , + \infty, \\
\nu_i (P_{k(f_i, f_j)} ^{(i)} ) & \mbox{ if } &  P_\l ^{(j)} \notin \theta_{fh} (f_i),\\
0 & \mbox{ if }&  \l = 0,\\
+ \infty & \mbox{ if } & j= i \mbox{ and }  \l = + \infty.
\end{array}
\right.
\end{equation}
\end{Defi}

\begin{Rem} \label{value}
If $h \in  {\cal RC}_f$ is irreducible then 
$\nu_i (P_{k(h,f)}^ h ) = \frac{ \r (f_i, h)  } {\deg h}$.
\end{Rem}

If $h \in {\cal C}_f$ and in particular when $h = f_Y$, 
we study the  $f$-bunch decomposition of $h$ 
by analysing the Newton polyhedra with respect to suitable coordinates, 
of the polynomials defining the images of $h= 0$ under 
 {\it quasi-ordinary morphisms} associated to the irreducible factors of $f$.

If $f \in \C \{ X_1, \dots, X_d \} [Y]$ is a quasi-ordinary polynomial we say that the morphism
\[
\left\{
\begin{array}{l}
\xi_f \;\; : \;\; (\C^{d+1}, 0)  \;\; \longrightarrow \;\;  (\C^{d+1}, 0)
\\
U_1 = X_1, \; \; \dots , \;\; U_d = X_d,  \; \; T = f(X_1, \dots,X_d, Y)
\end{array}
\right.
\]
is {\it quasi-ordinary}. By definition, {\it the critical space} of the morphism $\xi_f$ is
the {\it polar hypersurface}, $f_Y = 0$, associated to the 
given quasi-ordinary projection 
$(X_1, \dots, X_d, Y) \mapsto (X_1, \dots, X_d)$. 
The {\it discriminant space} is the image of the the critical space by $\xi_f$, see \cite{TD}.

More generally, if $h \in \C \{ X \} [Y]$
the equation defining the 
image of the hypersurface $h= 0$ by $\xi_f$
is obtained by eliminating $X_1, \dots, X_d, Y$  from the equations:
$h= 0$, $T- f= 0$,$U_1 = X_1, \dots, U_d =X_d$,
i.e., by the vanishing of:
\begin{equation}
\psi_f (h)  := \mbox{Res}_Y (T-f, h).
\end{equation}  
The degree of the polynomial $\psi_f ({h})  \in \C \{ U \} [T] $
is  equal to $\deg h$. 
If $h=0$ is analytically irreducible at the origin
the same holds for its image $\psi_f (h) =0$ thus $\psi_f (h)$ is an 
irreducible polynomial.
If $h = h_1 \cdots h_t$ then it follows from (\ref{Res}) that   $\psi_f (h) = 
\prod_{r=1}^t \psi_f (h_r)$.

We analyse the {Newton polyhedron}   of 
$\psi_{f_i} (h)$ for $f_i$ any irreducible factor of $f$.
Recall that the Newton polyhedron    ${\cal N} (\f) \subset \R^d$ 
of a non zero series 
$\f = \sum c_\a X^{\a} \in \C\{ X \} $ with $X=(X_1, \dots, X_{d})$
is the convex hull of the set 
$\bigcup_{c_\a \ne 0} \  \a  + \R_{\geq 0}^{d}$.
The Newton polyhedron of a polynomial  $F \in \C \{ X \} [Y ]$ is 
the polyhedron ${\cal N}(F) \subset \R^{d} \times \R$ of $F$ viewed as a series in 
$X_1, \dots, X_d, Y$. We need the following notation:

\begin{Not}
We denote the 
Newton polyhedron of $Y^p - X^a \in \C \{ X_1, \dots, X_d \} [Y]$
by the symbol 
$\frac{p}{\overline{q}}$ where $q:= \frac{a}{p} \in \Q^d$ is the 
{\it inclination} of the edge of ${\cal N} ( Y^p - X^a )$.
Our notation $\frac{p}{\overline{q}}$ is inspired by 
the one $\{ \frac{b}{\overline{a}} \}$  used by Teissier 
with a different meaning ($\{ \frac{b}{\overline{a}} \} := \frac{b}{\overline{\frac{a}{b}}}$)
to describe elementary Newton polygons in \cite{TD} and \cite{TJ}.
We have the following property of the Minkoski sum: 
$\frac{\underline{p}}{{q}} +   \frac{\underline{p'}}{{q}}
= \frac{\underline{p+ p'}}{{q}}$. 
\end{Not}

We prove that   the tree 
$\theta_f (f)$ determines 
the Newton polyhedra ${\cal N} (\psi_{f_i} (f_Y))$ and ${\cal N} (\psi_{f} (f_Y))$.
\begin{The} \label{New}   Let $f \in \C \{ X_1, \dots, X_d \} [Y]$  
be a quasi-ordinary polynomial with irreducible factors $f_1, \dots, f_s$.       
The Newton polyhedron of $\psi_{f_i} (f_Y) $ (resp. of $\psi_{f} (f_Y) $) is the Minkowski sum:
\begin{equation} \label{New2}
{\cal N} (\psi_{f_i} (f_Y)) =  \sum_j  \frac {c_j}{\overline{\nu_i (P_j)}}
\quad  \quad
\quad  \quad  \left( \mbox{ resp. } {\cal N} ( \psi_{f} (f_Y)) =  \sum_j \frac{ c_j }{\overline{\nu_1 (P_j) + \cdots +  \nu_s (P_j)}} \, \right),\end{equation} 
where in both cases $P_j$ runs through the set of non extremal vertices of $\theta_f (f)$.
\end{The}

We use this result to characterize the type of the bunch decomposition of 
the polar hypersurface $f_Y = 0$ in terms of the tree $\theta_f(f)$ and conversely, 
generalizing a theorem of  Garc\' \i a Barroso's 
for a generic polar curve of a  plane curve germ (see Th\'eor\`eme 6.1 of  \cite{GB}).
\begin{The} \label{main}
Let $f \in \C \{ X_1, \dots, X_d \} [Y]$  
be a quasi-ordinary polynomial with irreducible factors $f_1, \dots, f_s$.
\begin{enumerate}
\item
The type of the bunch decomposition of the partial derivative $f_Y$ is:
\begin{equation} \label{ftype}
\{ (\nu_1 (P_j) , \dots, \nu_s (P_j) ; \, \, c_j ) \}_j 
\end{equation}
where $P_j$ runs through the set of non extremal vertices of $\theta_f (f)$
and $c_j$ is the coefficient of $P_j$ in the chain $- \partial \g_f$.
In particular, when $f$ is irreducible with characteristic exponents 
$\l_1, \dots, \l_g$ the type of $f_Y$ is 
$$ \{   (\nu( P_{\l_j} ) ; \, \, n_0 n_1 \cdots n_{j-1}(n_j -1)) \}_{j=1}^g. $$
\item \label{gg}
The type of the bunch decomposition of $f_Y$ and the degrees of the irreducible factors of $f$ 
determine the Eggers-Wall tree of $f$. 
\end{enumerate}
\end{The}

\begin{Rem} \label{ppp} 
Assertion 1 of Theorem \ref{main} generalizes Popescu-Pampu's Theorem 3.8.5 of \cite{PP2} 
(or Theorem 6.3 of \cite{PP1})
 obtained in the case of a 
quasi-ordinary derivative (when $f_Y \in {\cal RC}_f$).
\end{Rem} 
Popescu Pampu's theorem  is based on a generalization of a result of Kuo and Lu \cite{KL}  Lemma 3.3,
which is also essential to prove the properties of the bunch decomposition of the polar curve given by
Eggers, Garc\'\i a Barroso and Wall (see \cite{Eggers}, \cite{GB} and \cite{Wall}).
Kuo-Lu's lemma compares 
in the case of a plane curve germ,  $F(X, Y) =0$,  
the dominant terms of the differences of any fixed root $Y =\z(X)$  of 
$F$ 
with the other roots of $F$ and of $\z(X)$  with the roots of $F_Y$.
The additional hypotesis needed to generalize Kuo-Lu's lemma to the case of a
quasi-ordinary polynomial $f\in \C \{ X_1, \dots, X_d \}[Y]$,
 i.e., to compare the 
{\it roots} of the derivative $f_Y$
with the roots of $f$,  is that the polynomial $f_Y$ should be 
radically comparable to $f$. 
With this hypothesis Popescu-Pampu's decomposition
follows 
by extending in a natural way the approach of Wall in the plane curve case
(see \cite{Wall}).
If $f_Y \notin  {\cal RC}_f$ the quasi-ordinary projection 
$(X_1, \dots, X_d, Y) \mapsto (X_1, \dots, X_d)$ may be replaced 
by a base change defined by an embedded resolution of 
the discriminant $\D_Y (f_Y\cdot f) = 0$, in such a way that the transforms of $f$ and $f_Y$ become
simultaneously quasi-ordinary with respect to the same quasi-ordinary projection over 
any point of the exceptional divisor. However, it is not clear that the decompositions obtained in this 
way come from a decomposition of $f_Y = 0$ since base changes
do not preserve irreducible components in general.

\section{Newton polyhedra  and toric geometry} 

We introduce in the following subsections the tools needed to prove 
the main results.

\subsection{Polygonal Newton polyhedra and their dual Newton diagrams}

If $\f \in \C\{ X \} $ is a non zero series 
in the variables $X = (X_1, \dots, X_d)$ we have that 
any linear form  $w \in (\R^{d})^*$ in the cone $ \D_{d} :=(\R^{d})^*_{\geq 0}$
defines a {\it face} ${\cal F}_w$  of the  polyhedron ${\cal N} (\f)$:
$$
{\cal F}_w = \{  v \in {\cal N} (\f) \, /  \, 
\langle w, v \rangle   = \inf_{v' \in {\cal N} (\f)} \langle w, v' \rangle \}.$$
All faces of the polyhedron ${\cal N} (\f)$ can be recovered in this way.
The face of ${\cal N} (\f)$  
defined by $ w$ is compact if and only if  $w $ belongs to the interior
 $\stackrel{\circ}{\D}_{d}$ of the cone 
$\D_{d}$.
The {\it cone $\s( {\cal F} ) \subset \D_d$  associated to the face $\cal F$}
of the polyhedron ${\cal N} (\f)$ is 
$\s( {\cal F} ) := \{ u \in \D_d  \; /
\; \forall v \in {\cal F}, 
\langle u, v \rangle   = \inf_{v' \in {\cal N} (\f)} \langle u, v'
\rangle \}$.
The {\it dual Newton diagram} $\Sigma ({\cal N} (\f))$
is the set of cones $\s( {\cal F})$, for $\cal F$ running through the set
of faces of the polyhedron ${\cal N} (\f)$ (see \cite{Kho}).
\begin{Rem} \label{lity}
If $\f = \f_1 \dots \f_r$,   the elements of the dual Newton diagram of $\f$ 
are the intersections $\cap_{i=1}^r \s_i$ for $\s_i$ running through $\Sigma ({\cal N} (\f_i))$ 
for $i =1, \dots, r$.
\end{Rem}
We deduce this property by duality from: 
\begin{equation} \label{fan}
{\cal N} (\f) = {\cal N} (\f_1) +  \cdots + {\cal N} (\f_r).
\end{equation}

The set of compact faces of a polygonal Newton polyhedron ${\cal N} (\f)$ is 
combinatorially 
isomorphic to a finite subdivision of a compact segment: since 
${\cal N} (\f)$ is polygonal 
the cones of the dual Newton diagram which intersect 
 $\stackrel{\circ}{\D}_{d}$
are of dimensions $d$ and $d-1$ by duality.

\begin{Lem} \label{cri}
If  $\f \in \C \{ X \}$ has a polygonal Newton polyhedron 
any irreducible factor of $\f$ which is not associated 
to $X_i$, for $i=1, \dots, d$,
has a polygonal Newton polyhedron.   
\end{Lem}
{\em Proof.}
It follows from (\ref{fan}) that the compact face of ${\cal N} (\f)$
determined by $w \in \stackrel{\circ}{\D}_d$ is 
the Minkowski sum of the compact faces,  determined by $w$, 
on the Newton polyhedra of the factors.
Since the polyhedron ${\cal N} (\f)$ is polygonal the dimension of these 
compact faces is zero or one. It follows that 
the Newton polyhedron of an irreducible factor of $\f$ is polygonal or
a translation of $\R_{\geq 0}^d$, and in the latter case
this irreducible factor is associated to one variable. 
\hfill $\ {\diamond}$

\subsection{Coherent polygonal paths}

The  Newton polyhedron  of a polynomial $0 \ne F \in \C \{ X \} [Y ]$
is contained in $ \R^{d} \times \R$.
Any  {\it irrational} vector $w \in \D_d$, i.e.,  
with  linearly independent coordinates over $\Q$, defines a {\it  coherent polygonal path}
on the compact edges of ${\cal N} (H) $ (the terminology comes from the combinatorial convexity theory,
 see \cite{Billera}).
This path is defined by $c_w (t) = (u_w(t), t)$ 
where $u_w (t) $ is the unique  point of the hyperplane section $v= t$ of 
${\cal N} (H) $ where  the minimal value of the linear function $w$ is reached
for $t \in [\mbox{ord}_Y H, \deg H]$ (for $\mbox{ord}_Y $ the order of $F$ as a series in $Y$).
The point $u_w (t)$ is unique because the vertices of the polyhedron ${\cal N} (H) $ 
are rational, i.e., they belong to the lattice $\Z^d \times \Z$.
Any compact edge which is not parallel to the hyperplane $v= 0$ belongs to
some path $c_w (t)$ for some irrational vector $w \in \D_d$.
The maximal segments of the polygonal path $c_w (t) $ 
are of the form 
$\e_i = [p_i , p_{i+1}] $ where $p_j = (u_j, v_j)$ for $j= i, i+1$ and $v_i < v_{i+1}$.
We call the vector $q_i = \frac{u_i - u_{i+1}}{v_{i+1} - v_i}$ the {\it inclination}
and the integer $l_i = {v_{i+1} - v_i}$ the 
{\it height} of the edge $\e_i$ (see \cite{GP1}  where this construction 
is related to generalizations of Newton Puiseux Theorem).

\begin{Lem} \label{plo} 
Let  $\{ u_i \}_{i=1}^r$ be $r$ different non zero vectors in  $\Q^d$
such that 
$0 < u_r \leq \cdots \leq u_1$ (with respect to the order (\ref{order}))
and integers $l_1, \dots, l_r \in \Z_{> 0}$.
The Minkowski sum: 
\begin{equation} \label{Min}
{\cal N} = \sum_{i=1}^{r} \frac{l_i}{\overline{u_i}}
\end{equation}        
is a polygonal polyhedron in $\R^{d+1}$.
It has $r$ compact edges ${\cal E}_i$ of inclinations $u_i$ and 
heights $l_i$ for $i=1, \dots, r$. 
The polyhedron ${\cal N}$ determines the terms of the 
Minkowski sum (\ref{Min}).  
\end{Lem}
{\em Proof.}
The vector hyperplane $h_i$ orthogonal to the compact edge of $\frac{l_i}{\overline{u_i}} $
defines two half-spaces which subdivide the  interior of the cone 
$\D_{d+1}$ since $0 < u_i$. The condition $u_i < u_j$ implies that 
the hyperplanes $h_i$ and $ h_j$ 
do not intersect in the interior $\stackrel{\circ}{\D}_{d+1}$ of the cone 
$\D_{d+1}$. It follows that the possible codimensions of the cones of
the dual diagram of ${\cal N}$, intersecting $\stackrel{\circ}{\D}_{d+1}$, 
are $0$ and $1$. The codimension one case corresponds to the cones 
defined by the hyperplane sections  $h_i$.
By duality the polyhedron ${\cal N}$ is polygonal. The edge defined 
by $u \in h_i \cap \stackrel{\circ}{\D}_{d+1} $ is the Minkowski sum
of the faces defined by $u$ on each of the terms of (\ref{Min}), i.e, 
it is a translation of the polyhedron $\frac{l_i}{\overline{u_i}}$. 
\hfill $\ {\diamond}$

\begin{Lem} \label{curve}
Let $F \in \C \{ X \} [Y] $ be a monic polynomial of degree $>0$ with $0 \ne F(0)$ a non unit.
If the polygonal path $c_w (t)$ does not depend on the irrational vector of $\w \in \D_d$
then the inclinations of the edges of $c_w (t)$ are 
totally ordered with respect to the order (\ref{order}) and 
the polyhedron ${\cal N} (F)$ is polygonal.
\end{Lem}
{\em Proof.} 
We label the edges of $c_w (t)$ by $\e_0, \dots, \e_r$ 
in such a way that $v_0 = 0 < v_1 < \cdots < v_{r} < v_{r+1} = \deg F$ with the previous notations.
The irrational vector $w$ defines the total order of $\Q^d$ defined by 
$u \leq_w u' \Leftrightarrow \langle u, w \rangle  \leq  \langle u', w \rangle$
and we have that 
\begin{equation}\label{w}
q_r <_w q_{r-1} <_w \cdots <_w q_0
\end{equation}
(see Lemme 5 of \cite{GP1}). 
By hypothesis the path $c_w (t)$ 
does not depend on the irrational $w \in \D_d$. 
It follows that the inequality  (\ref{w}) holds for all irrational vector $w \in \D_d$
therefore  $q_r \leq q_{r-1} \leq \cdots \leq  q_0$ with respect to the order (\ref{order}).
It follows that the polyhedron ${\cal N} (F)$ is of the form 
(\ref{Min}) therefore it is polygonal by lemma \ref{plo}. 
\hfill $\ {\diamond}$

\subsection{A reminder of toric geometry}
We give some definitions and notations (see  
\cite{Ewald}, \cite{Oda} or \cite{TE} for proofs).
If $N \cong \Z^{d}$ is a lattice we denote by $M$ the dual lattice, by
$N_\R$ the real vector space spanned by $N$. 
In what follows a {\it cone} mean a  
{\it rational convex polyhedral cone}:
the set of   non
negative linear combinations of vectors $a^1, \dots, a^s  \in N$.
The cone $\sigma$ is {\it strictly convex} if
$\sigma$ contains no linear subspace of dimension $>0$;
the
cone $\sigma$ is {\it regular\index{regular cone}} if
the primitive integral vectors defining the $1$-dimensional faces
belong to a basis of the lattice  $N$.
The {\it dual} cone  $\sigma^\vee$ (resp. {\it orthogonal} cone 
$\sigma^\bot$) of $\sigma$ is the set
$ \{ w  \in M_\R / \langle w, u \rangle \geq 0,$  (resp. $ \langle w, 
u \rangle = 0$)  $ \; \forall u \in \sigma \}$.
A {\it fan\index{fan}} $\Sigma$ is a family of strictly convex
  cones  in $N_\R$
such that any face of such a cone is in the family and
the intersection of any two of them is a face of each.
The {\it support\index{support of a fan}} of the fan $\Sigma$ is the set
$ \bigcup_{\sigma \in \Sigma} \sigma \subset  N_\R$.
The fan $\Sigma$ is {\it regular} if all its cones are regular.
If  $\s$ is a cone in the fan $\Sigma$ 
the semigroup $\s^\vee  \cap M$  is of finite type,  it spans the lattice $M$ and 
defines the affine variety 
$Z^{\s^\vee \cap M} = \makebox{Spec} \, \C [\s^\vee \cap M  ]$, which 
we denote also by  $Z_{\s, N}$ or by $Z_\s$ when the lattice is clear from the context.
If $\s \subset \s'$ are cones in the fan $\Sigma$ then
we have an open immersion $Z_\s \subset Z_{\s'} $;  
the affine varieties $Z_\s$  corresponding to cones in a fan
$\Sigma$ 
glue up to define the {\it toric variety\index{toric variety}} 
$ Z_\Sigma$. The toric variety $Z_\Sigma$ is non singular if and only
if the fan $\Sigma$ is regular.
The torus,   $(\C^*)^{d+1} $,   
is embedded as an open dense subset $Z_{ \{ 0 \} }$ of  
$ Z_\Sigma$, which acts  on each chart $Z_\s$; these actions
paste to an action on $ Z_\Sigma$  which extends
the action of the torus on itself by multiplication.
The correspondence which associates to a cone  $\s \in \Sigma$ 
the Zariski closed subset  $\O_{\s}$ of  $Z_\s$, 
defined by the ideal $(X^w/ w  \in (\s^\vee - \s^\bot) \cap M
)$ of $\C[\s^\vee  \cap M ]$,
is a bijection between $\Sigma$ and the set of orbits of the torus action in $Z_\Sigma$.
For example the set of faces of a cone $\s$ defines a fan such that the associated toric 
variety coincides with $Z_\s$.

We say that a fan $\Sigma'$ is a {\it subdivision\index{fan subdivision}} of the fan
$\Sigma$ 
if both fans have the same support and if every cone of
$\Sigma'$ is contained in a cone of $\Sigma$.
If $ \Sigma' \ni \s' \subset \s \in  \Sigma$  we have the morphism 
$Z_{\s'} \rightarrow Z_\s$ defined by 
the inclusion of semigroups $\s^\vee \cap M  \rightarrow {\s'}^\vee \cap M$.
These morphisms glue up and define the {\it toric
modification\index{modification}}
$ \pi_{\Sigma'} : Z_{\Sigma '} \rightarrow   Z_\Sigma$.
Given any fan $\Sigma$  there exists a regular fan $\Sigma'$ subdividing $\Sigma$ (see \cite{Oda}). 
The associated toric modification $ \pi_{\Sigma'}$ is a desingularization.

For instance, if we denote the lattice $\Z^{d}$ by $N$ 
we have that   $\D_d \subset N_\R$,
the toric variety $Z_\D$ is the affine space $\C^d$ and the orbits 
correspond to the strata of the stratification defined by the 
coordinate hyperplanes $X_i = X^{e_i} = 0$
(where $e_i$ runs through the elements of the basis of $M$ defined by the edges of $\D_d^\vee$).
Any fan $\Sigma$ supported on $\D_d $   
defines the toric modification $\p_{\Sigma} : Z_\Sigma \rightarrow \C^{d}$.
Taking away  the cone $\s$ from the fan of the cone $\s$ means geometrically to
take away the orbit $\O_\s$ from the variety $Z_\s$.
This implies that the exceptional fiber  $\p_\Sigma^{-1} (0) $
of the toric modification  $\p_\Sigma$ is 
\begin{equation} \label{exceptional}
\p_\Sigma^{-1} (0) = \p_{\Sigma}^{-1}( \O_{\D_d} ) = 
\displaystyle{\bigcup_{\t \in \Sigma, \stackrel{\circ}{\t} \subset \stackrel{\circ}{\D_d }}} 
\O_\t
\end{equation}
Let ${\cal V}$ be a subvariety of $\C^d $ such that the intersection with 
the torus is a non singular dense open subset of ${\cal V}$. 
The {\it strict transform} ${\cal V}_\Sigma \subset Z_\Sigma $  
is the subvariety of $\pi_\Sigma^{-1} ({\cal V})$ 
such that the restriction $ {\cal V}_\Sigma  \rightarrow {\cal V}$ is a modification. 
If  $0 \ne \f = \sum c_a X^a \in  \C \{ X \}$ is a non zero series in $X = (X_1, \dots, X_d)$
the dual Newton diagram $\Sigma( {\cal N} (\f))$ is a subdivision of 
$\D_d$.  
The {\it symbolic restriction} of $\f$ to a set ${\cal F} \subset M_\R$ is 
$\f_{| \cal F} :=\sum_{a \in {\cal F}} c_a X^{a} $.

\begin{Lem} \label{nondeg}
Let $\f = \sum c_a X^a \in  \C \{ X_1, \dots, X_d \}$ be an irreducible series, not associated to any 
of the variables $X_1, \dots, X_d$, defining the germ $({\cal V}, 0) \subset (\C^d , 0)$.
Let $\Sigma$ be any subdivision of the dual Newton diagram 
$ \Sigma ({\cal N} (\f))$.
If  $\s \in \Sigma$ and if 
$\stackrel{\circ}{\s} \subset \stackrel{\circ}{\D_d}$, 
the intersection 
$\O_\s \cap {\cal V}_\Sigma $ is defined by the vanishing of
$X^{-u} \f_{{\cal F}_\s} \in \C [ \s^\bot \cap M ]$
where $u$ is any vertex of the compact face ${\cal F}_\s$ 
of ${\cal N} (\f)$ defined by 
$\s$.
\end{Lem} 
{\em Proof.} Let $v \in M $ such that $-v + {\cal N} (\f) \subset \s^\vee$. 
Then all the terms in $X^{-v} \f$ vanish on the orbit $\O_\s$ (since their exponents belong to $\s^\vee - \s^\bot$)
unless the vector $v$ belongs to  the affine hull $\mbox{Aff} ( {\cal F}_\s ) $ of the compact face ${\cal F}_\s$. 
In this case we have that $(X^{-v} \f) |{\O_\s}  = X^{-v} \f_{| {\cal F}_\s}$.
If $v, v' \in \mbox{Aff} ( {\cal F}_\s ) \cap M$ we have that $v- v' $ belongs to
$\s^\bot \cap M$ therefore the polynomials  $X^{-v} \f_{| {\cal F}_\s}$ and $X^{-v'} \f_{| {\cal F}_\s}$
are related by the invertible function $X^{v-v'}$ on the torus $\O_\s$.
It follows from this that $X^{-v} \f_{| {\cal F}_\s}$ defines the ideal of the intersection
${\cal V}_\Sigma \cap \O_\s$.
\hfill $\ {\diamond}$

\subsection{An irreducibility criterion for series with
polygonal Newton polyhedra} \label{ch}

We use Theorem \ref{New} to prove Theorem \ref{main}
by translating the existence of
the bunch decomposition of $f_Y$ in geometrical terms
by means of an {\it irreducibility criterion} for power series
with {\it polygonal} Newton polyhedron:

\begin{Defi} \label{like}
A polyhedron 
is polygonal if the maximal dimension of 
its compact faces is one.
\end{Defi}
Polygonal polyhedra share some properties with classical Newton polygons of
plane curve germs. For instance, any Newton polygon is the Minkowski sum of 
{\it elementary Newton polygons} up to translations (see \cite{TJ}).
The criterion holds when the field $\C$ of 
complex numbers is replaced by an algebraically closed field 
of arbitrary characteristic.
The criterion generalizes a fundamental property of plane curves.
We introduce some definitions and notations:

Let  ${\cal E} = [\a, \a']$ be a compact segment joining two elements $ \a, \a'$ of the lattice 
$\Z^d$, and denote by $u$ the primitive integral vector parallel to $ {\cal E}$, i.e., 
we have an equality of the form  $\a' - \a = l u$ for some maximal integer $l \in \Z_{\geq 1}$.
If $\f = \sum c_a X^a $ is a series in $\C \{ X_1, \dots, X_d \}$ we have that 
\begin{equation} \label{edge}
X^{-\a} ( \sum_{a \in {\cal E} \cap \Z^d} c_a X^a)  
= \sum_{i=0}^{l} c_{iq} X^{i q} = p({\cal E}, \f)(X^u)
\end{equation}
where $p(\f,{\cal E})$ is the polynomial $p(\f,{\cal E}) = \sum_{i=0}^{l} c_{iq} t^i$.
Obviously, this definition depends on the 
order of the vertices: the polynomial obtained by interchanging 
the vertices $\a$ and $\a'$ of ${\cal E}$ 
is equal to  $t^l p(\f, {\cal E}) (t^{-1})$.
In both cases  
these two polynomials  
define isomorphic subschemas of $\C^*$. 
If in addition $\f $ is a monic polynomial in $\C \{ X_1, \dots, X_{d-1} \} [X_d]$
we fix the order of the vertices  $\a= (a, t)$ and  $(a', t') $ of  ${\cal E}$
by the convention $t < t'$, and we obtain that the polynomial $p(\f,{\cal E})$ is 
defined in a unique way by (\ref{edge}).

\begin{The}\label{irred}
If $\f \in \C \{ X_1, \dots, X_d \}$ is irreducible and has a polygonal Newton polyhedron ${\cal N} (\f)$, 
then the polyhedron  ${\cal N} (\f)$ 
has only one compact edge ${\cal E}$ and 
the polynomial $p(\f,{\cal E})$ has only one root in $\C^*$.
\end{The} 
{\it Proof.} 
Since ${\cal N} (\f)$ is polygonal the cones  $\s$ in the  
dual Newton diagram $\Sigma$ of the polyhedron  ${\cal N} (\f))$, 
such that $\stackrel{\circ}{\s} \subset \stackrel{\circ}{\D_d}$ 
are of codimensions  $0$ or $1$ (the possible dimensions
of the compact faces of ${\cal N} (\f)$).
We keep notations of lemma   \ref{nondeg}.
If $\dim \s = d$ it follows from lemma  \ref{nondeg} 
that the intersection ${\cal V}_\Sigma \cap \O_\s$
is empty.
If $\dim \s = d-1$ the cone $\s$ corresponds to the compact edge $\cal E$ of the polyhedron ${\cal N} (\f)$.
By lemma  \ref{nondeg} 
the intersection ${\cal V}_\Sigma \cap \O_\s$ is defined by the vanishing of  
$X^{-u}\f_{\cal E}$ on the torus $\O_\s$ (the vector $u$ being one of the vertices of the edge $\cal E$).
We have that the coordinate ring $\C [\s^\bot \cap M]$  of the orbit $\O_\s \cong \C^* $ 
is isomorphic to $\C [ X^{\pm u}]$, where $u$ the primitive integral vector parallel to 
the edge  ${\cal E}$.  
By formula 
(\ref{edge}) the polynomial $X^{-u}\f_{\cal E}$ corresponds to the polynomial 
$ p({\cal E}, \f)(X^u) $. It follows that the intersection 
${\cal V}_\Sigma \cap \O_\s$ is a finite set of points counted with multiplicities
which correspond to the zeroes of the polynomial $ p({\cal E}, \f)$.
By (\ref{exceptional})  the fiber of the modification 
$\p_{\Sigma}| {\cal V}_\Sigma : {\cal V}_\Sigma \rightarrow {\cal V}$
is equal to the discrete set   $\bigcup (\O_\s \cap {\cal V}_\Sigma)$, where $\s$
runs through the cones 
$\s \in \Sigma$ such that $\stackrel{\circ}{\s} 
\subset \stackrel{\circ}{\D_d}$.
Since by hypothesis the germ ${\cal V}$ is analytically irreducible at the origin 
this fiber is a connected set by Zariski's Main Theorem 
(see \cite{Mumford} and \cite{Main}) thus it is reduced to one point. 
This implies that the Newton polyhedron ${\cal N} (\f)$ has 
only one compact edge ${\cal E}$ and 
that the polynomial $p({\cal E}, \f)$ has only one root in $\C^*$.  
\hfill $\ {\diamond}$

We will need the following lemma in section \ref{le}.
We denote by $M$ (resp. by $M'$) the lattice spanned by the exponents of 
monomials in $\C \{ X_1, \dots, X_d \} $ (resp. in $\C \{ X_1, \dots, X_{d}, Y \} $).
\begin{Lem} \label{spe}
Let $h,  h' \in \C \{ X_1, \dots, X_{d} \} [Y]$ be monic polynomials.
We suppose that the polyhedron ${\cal N} (h')$ is polygonal with compact edge of 
the inclination $\l \in M_\Q$. We
denote by $\cal E$ the compact edge of the polyhedron $\frac{ \underline{ \deg h }}{\l} $.
If 
${\cal N} (h) \subset \frac{ \underline{ \deg h }}{\l} $ and if $p( {\cal E} , h) = t^{\deg h} $, 
the strict transform of $h = 0$ by $ \p_ \Sigma$, for $\Sigma = \Sigma (\frac{ \underline{ \deg h }}{\l}   )$, 
intersects $ \p_\Sigma ^{-1} ( 0) $ only at the zero dimensional orbit $\O_\t$ 
where $\t$ is the cone associated to the vertex $ ((\deg h)  \l , 0)$ of the polyhedron
$\frac{ \underline{ \deg h }}{\l}$.
\end{Lem}
{\em Proof.} 
Denote by $v$  the vertex $ (0, \deg h) $ of $\frac{ \underline{ \deg h }}{\l}$.
We deduce from the hypothesis ${\cal N} (h) \subset \frac{ \underline{ \deg h }}{\l} $ and  $p( {\cal E} , h) = t^{\deg h} $ that
the series $Y^{-\deg h} h$ 
has terms in $ \C [\s^\vee \cap M' ] $,  constant term 
equal to one and the exponents of non constant terms do not belong
to $\s^\bot$,  
for 
$\s = \s ( \{ v \} )$ or $\s = \s ( {\cal E} )$.
It follows that the strict transform of $h= 0$ does not meet $\O_\s$ since 
the terms $\ne 1$  
vanish on $\O_\s$. 
We deduce from (\ref{exceptional}) that the strict transform of $h$ intersects 
the exceptional fiber $ \p_\Sigma ^{-1} ( 0) $ at the closed orbit 
$\O_\t$  (which is reduced to a point since $\t$ is of maximal dimension $d+1$).
\hfill $\ {\diamond}$

\section{The proofs of the results on the type of $f_Y$ \label{tbc}}

We apply the irreducibility criterion to clarify the relation between
the type of the $f$-bunch decomposition of $h \in {\cal C}_f$ and 
the Newton polyhedra ${\cal N} (\psi_{f_i} (h))$.

\begin{Lem} \label{monotone} 
The restriction of the valuation $\nu_i$ to $\theta_{fh} (f_i) $ is an order-preserving
bijection. The characteristic exponents of $f_i$ and 
the  valuation $\nu_i (P_\l^{(i)})$ determine $\l$.
\end{Lem}
{\em Proof.}
We denote the characteristic exponents of $f_i$ by  $\l_1, \dots,\l_g$.
 Let  $P_\l $ be a vertex of $ \theta_{fh} (f_i)$.  
For simplicity we drop the index $i$.
If $P_\l <  P_{\l_1}$  then we have that $\nu (P_\l) = n \l $.
Otherwise there exists a unique $1 \leq j \leq g $ 
such that $P_{\l_j} \leq P_\l < P_{\l_{j+1}}$ since
$ \theta_{fh} (f_i)$ is totally ordered.
Then the first assertion follows from the inequality: 
\[
\begin{array}{lll}
\nu (P_{\l_j}) &  = & e_{j-1} \l_j  + \sum_{k=1}^{j-1} (e_{k-1} - e_{k}) \l_k   
  =   e_{j} \l_j  + \sum_{k=1}^{j} (e_{k-1} - e_{k}) \l_k   \leq \\
& \leq & e_{j} \l + \sum_{k=1}^{j} (e_{k-1} - e_{k}) \l_k    
=  \nu (P_{\l})  <  e_{j} \l_{j+1} + \sum_{k=1}^{j} (e_{k-1} - e_{k}) \l_k =  
\nu (P_{\l_{j+1}}).
\end{array}
\]
If we know the characteristic exponents of $f_i$
then there is a unique $j$ such that $ \nu (P_{\l_j}) \leq \nu (P_\l) < \nu (P_{\l_{j+1}})$
where  we convey that $\l_{0} =0$. Then we recover $\l$ from  equation (\ref{nu}).
\hfill $\ {\diamond}$

If $h \in {\cal RC}_f$ we show below that the Newton polyhedron of 
$\psi_{f_i} (h)$ is determined by the $f$-type of $h$.
\begin{Pro} \label{Newton1}
If $\{ (q_{1,r}, \dots, q_{s,r} ; c_r) \}_{r=1}^{s(f,h)}$ is the type of the
$f$-bunch decomposition of a polynomial 
$h \in {\cal RC}_f$
then the Newton polyhedron of $\psi_{f_i} (h) $ is
the Minkowski sum:
\begin{equation} \label{splut}
{\cal N} (\psi_{f_i} (h)) =  \sum_{r = 1}^{s(f,h)} \frac{c_r}{\overline{q_{i,r} }}.
\end{equation}
\end{Pro}
{\em Proof.}
If $\{ \t_r^{(j)} \}_{j=1, \dots, \deg b_r}$ are the roots of the factor $b_r$ of the 
$f$-bunch decomposition of $h$ for $r=1, \dots, s(f,h)$, viewed in 
some suitable ring extension of the form $\C \{ X^{1/k} \}$,  
it follows from proposition \ref{coincidence}
and the definition of the bunches that:
\begin{equation} \label{splot}
f_i ( \t_r^{(j)} ) = X^{ q_{i,r} } \epsilon_{i, r, j},   \mbox{ where  }  \epsilon_{i, r, j} 
\mbox{ is a unit in   } \C \{ X^{1/k} \}.
\end{equation}
By general properties of the resultant we have that:
$
\psi_f ({b}_r) = \prod_{j=1}^{\deg h} (T - f( \t^{(j)}_r))
$.
We deduce from this  and (\ref{splot}) that the polyhedron
${\cal N} ( \psi_{f_i} (b_r) ) = \sum_{j=1}^{c_r}   {\cal N} ( T - f_i ( \t_r^{(j)} ) )$
is equal to $\frac{c_r}{\overline{q_{i,r} }}$
and equality (\ref{splut}) follows from the property (\ref{fan}).
\hfill $\ {\diamond}$


\begin{Rem} \label{Newton3}
If  $h \in {\cal C}_f$ the assertion of proposition 
\ref{Newton1} is not true in general (see remark \ref{ple}).
\end{Rem}

The following proposition generalizes Proposition 3.4.8 \cite{PP2}.
\begin{Pro} \label{Wall} {\rm }
\noindent
Given $\theta_f (f)$, if $h \in {\cal RC}_f $ 
the following informations determine each other:
\begin{enumerate}
\item The contact chain $[h]^{(f)}$.
\item The type of the bunch decomposition of $h$ induced by $f$.
\item The collection of Newton polyhedra of $\psi_{f_i} (h)$ for $i=1, \dots,s$. 
\end{enumerate}
\end{Pro}
{\em Proof.} 
$1. \Rightarrow 2.$ : If the contact chain is  $[h]^{(f)} = \sum c_i P_i$
the type of the $f$-bunch decomposition of $h$ is 
$\{ ( \nu_1(P_i), \dots, \nu_s (P_i) ; c_i) \}_i$
by remark \ref{value}.
The implication $2. \Rightarrow  3.$ is a direct consequence of 
proposition \ref{Newton1}. 
$3. \Rightarrow  1.$ :
The Newton polyhedron of $\psi_{f_i} (h)$ is polygonal by proposition \ref{Newton1}.
We recover the set of  vertices  
$\{ P_j ^{(i)} \}_j $ of $\theta_{fh} (f_i)$ corresponding to the  
orders of coincidence of 
$f_i$ with the irreducible factors of $h$
from  
the inclinations of the compact edges  
by lemma \ref{monotone} (since $\theta_f (f_i)$ is given).
The maximal point  $P_{j_0} ^{(i)}$ of the set $\{ P_j ^{(i)} \}_j $
corresponds to a factor $b_1$ of the $f$-bunch decomposition
of degree $c_1$ equal to the  height of the edge of ${\cal N} (\psi_{f_{i}} (h)) $
of maximal inclination $\nu_{i} (P_{j_0} ^{(i)}) $.
Then we can replace $h $ by $h '= h/ b_1$ and continue in the same way.
The Newton polyhedra of $\psi_{f_i} (h')$ is obtained 
from ${\cal N} (\psi_{f_i} (h'))$ by subtracting the elementary polyhedra 
$ \frac{{c_1} }{\overline{\nu_{i} (P_{j_0} ^{(i)})}}$ (the subtraction makes sense
by lemma   \ref{plo}).
\hfill $\ {\diamond}$

\begin{Pro} \label{Wall2}
Given $\theta_f (f)$ and  $h \in {\cal C}_f $, if the Newton 
polyhedra $\psi_{f_i} (h)$ for $i=1, \dots,s$, are polygonal
the following informations determine each other:
\begin{enumerate}
\item The type of the bunch decomposition of $h$ induced by $f$.
\item The collection of Newton polyhedra of $\psi_{f_i} (h)$ for $i=1, \dots,s$. 
\end{enumerate}
\end{Pro}
{\em Proof.} 
If $\deg h = m > 0$ the type of $h$ determines two vertices of  ${\cal N} (\psi_{f_i} (h)) $:
the vertex $(0, m)$ which corresponds to the monomial $T^m$
and the vertex  $( \r(f_i, h), 0)$ which  corresponds to
the dominant term of $(\psi_{f_i} (h))_{|T = 0} = \mbox{Res}_Y (-f_i , h)$.
In particular if $h$ is irreducible and if  $\psi_{f_i} (h)$ has a polygonal Newton polyhedron
then by the irreducibility criterion 
${\cal N} (\psi_{f_i} (h)) $ has only one compact edge equal to $[(0, m),( \r(f_i, h), 0)] $,  hence
we obtain that 
${\cal N} (\psi_{f_i} (h))  = 
\frac{\deg h}{ \overline{\frac{\r(f_i, h)}{\deg h}}}$.
In the general case, if the monic polynomial  
$\psi_{f_i} (h)$ has a polygonal Newton polyhedron the same holds 
for its irreducible factors by lemma \ref{cri}.
It follows then that the type of $h$ determines the Newton polyhedron of $\psi_{f_i} (h)$
(see the proof of proposition  \ref{Newton1}).
Conversely, if we are given the tree $\theta_f (f)$ 
and  polygonal polyhedron ${\cal N} ( \psi_{f_i} (h))$
for $i =1, \dots, s$ we recover the type from the inclinations 
and heights of the compact edges of ${\cal N} ( \psi_{f_i} (h))$
following the method of proposition \ref{Wall}.
\hfill $\ {\diamond}$

We determine the structure of  the Newton polyhedron of 
$\psi_{f_i} ( f_Y)$ by using  adequate 
toric base changes, which reduce to the case
$f_Y \in {\cal RC}_f$, 
in such a way that we can recover 
the  
coherent polygonal paths on the Newton polyhedron  of $\psi_{f_i} ( f_Y)$.

\begin{Lem} \label{cl}
If $f_Y \in {\cal RC}_f$ then the Newton  polyhedra of $\psi_{f_i} ( f_Y)$ are polygonal, for $i=1, \dots, s$. 
\end{Lem}
{\em Proof.} 
If $f_Y \in {\cal RC}_f$ the type of $f_Y$ is given 
in terms of the Eggers-Wall tree
by formula (\ref{ftype}) by applying Popescu-Pampu decomposition (see remark \ref{ppp}).
Then we obtain the polyhedra $\psi_{f_i} ( f_Y)$ for $i=1, \dots, s$, from the type of $f_Y$ 
by formula (\ref{splut}). 
These polyhedra are polygonal by lemmas \ref{monotone} and  \ref{plo}.
\hfill $\ {\diamond}$

We describe first  the toric base changes, already used in \cite{GP1} and \cite{Tesis}, 
we use to reduce to the quasi-ordinary derivative case.
The ring of 
convergent (or formal) complex power series in $X= (X_1, \dots, X_d) $
can be denoted by $\C \{ \D_d^\vee \cap M \}$
where $M $ denotes the lattice $\Z^d$ and $\D_d$ denotes the cone $(\R^d)^*_{\geq 0}$.  
The advantage of this notation is 
that we can define ring homomorphisms by changing the lattice $M$ 
or the cone $\D_d$. 

Let  $\t \subset \D_d$ be 
a regular\footnote{The regularity of the cone $\t$ is not essential 
in what follows (see \cite{GP1} and  \cite{GP3}).} cone of dimension $d$ 
(see \cite{Ewald} for its existence).
The dual cone $\t^\vee $ is generated by a basis $a_1, \dots, a_d$ of the lattice $M$ 
and contains $\D_d^\vee$. We define from the semigroup inclusion  
$\D_d^\vee \cap M \hookrightarrow \t^\vee \cap M$ the local ring extension $\C \{ \D_d^\vee \cap M \} \hookrightarrow \C \{ \t^\vee \cap M \}$. 
The local ring $\C \{ \t^\vee \cap M \}$  is equal to 
$\C \{ V_1, \dots, V_d \}$ where $V_i = X^{a_i}$ for $i =1, \dots, d$. 
We denote by  $H^{(\t)}$ the image of a polynomial
$H \in \C   \{ \D_d^\vee \cap M \} [Y]$
in the ring  $\C \{ \t^\vee \cap M \}[Y] $. 

By definition of  Newton polyhedron we deduce that: 
\begin{Rem}
Let $\f \in  \{ \D^\vee \cap M \}$ and $H \in \C   \{ \D_d^\vee \cap M \} [Y]$ be non zero. 
We have that: 
\begin{equation} \label{tau}
{\cal N} (\f^{(\t)} ) =   {\cal N} (\f) + \t^\vee,  \quad \mbox{ and } \quad 
{\cal N} (H^{(\t)})  =  {\cal N } (H )  + (\t^\vee \times \{ 0 \}).  
\end{equation}
\end{Rem}
We deduce from (\ref{tau}) the following lemma:
\begin{Lem} \label{tau2}
If $\t \subset \D_d$ is a regular cone 
the coherent polygonal paths defined by 
an irrational vector  $w \in \t$ on the edges of the polyhedra 
${\cal N}(F) $ and ${\cal N} (F^ {(\t)})$ coincide. \hfill $\ {\diamond}$
\end{Lem}

The main argument of the following proposition is already used in Proposition 2.14 of \cite{Tesis}.
\begin{Lem} \label{soy} 
If $\t \subset \D_d$ is a regular cone and if $f \in \C \{ X \} [Y]$ is
a quasi-ordinary polynomial then $f^{(\t)}$ is a quasi-ordinary polynomial
and $\tilde{\theta}_f (f) = \tilde{\theta}_{f^{(\t)}} ( f^{(\t)} )$.
\end{Lem}
{\em Proof.} 
If $\z^{(i)} \in \C \{ \D_d^\vee \cap \frac{1}{k} M \}$ is a root of $f$ 
then $ (\z^{(i)})^{(\t)}  \in \C \{ \t^\vee \cap \frac{1}{k} M \}$ is a root of $f^{(\t)}$.
Extending the cone does not modify the support of the series nor the lattices 
spanned by the exponents.
It follows then from Lipman's  characterization of
roots of quasi-ordinary polynomials that  $\z^{(i)}$ and $ (\z^{(i)})^{(\t)}$ 
have characteristic exponents defined by the same elements of the lattice $\frac{1}{k} M$,  
and that if $f$ is irreducible the same holds for $f^{(\t)}$
(see Proposition 1.5 of \cite{Lipman1} or Proposition 1.3  \cite{Gau}).
Then the equality $\tilde{\theta}_f (f)  =  \tilde{\theta}_{f^{(\t)}} ( f^{(\t)} )$ follows
from the fact that the ring extension 
$\C \{ \D_d^\vee \cap \frac{1}{k} M \} \hookrightarrow 
\C \{ \t^\vee \cap \frac{1}{k} M \}$
sends monomials 
to monomials and units to units.
\hfill $\ {\diamond}$

\begin{Rem} \label{subtil}
If $\t \ne \D_d$ then ${\theta}_f (f) \ne {\theta}_{f^{(\t)}} ( f^{(\t)} )$
since the valuations of the vertices, denoted 
$v$ and $v^{(\t)}$, 
take on any non extremal vertex $P$ of $\tilde{\theta}_f (f) = \tilde{\theta}_{f^{(\t)}} ( f^{(\t)} )$
values equal to  the coordinates of 
$\tilde{v}(P) = \tilde{v}^{(\t)}(P)$ with respect to two fixed different basis of the 
lattice $\frac{1}{k} M$. 
\end{Rem}

\noindent
{\bf Proof of Theorem  \ref{New}.} We discuss first the case of the Newton polyhedron  of $\psi_{f_i} (f_Y)$.
If $f_Y \notin {\cal RC}_f$, let $\Sigma$ a regular subdivision of 
the dual Newton diagram of $\D_Y (f_Y) $.
Let  $\t \in \Sigma$ be a cone of dimension $d$.
It follows that 
$(f^{(\t)} )_Y = (f_Y)^{(\t)}$ is polynomial in ${\cal RC}_{f^{(\t)} }$
by lemma \ref{soy} and the definitions.
By remark \ref{tau2}, the coherent polygonal path $c_w (t)$ 
determined by an irrational vector $w \in \t$ on 
${\cal N} (\psi_{f_i } (f_Y))$
and on ${\cal N} (\psi_{f_i^{(\t)}} (f_Y^{(\t)} ))$
coincide. 
By lemma \ref{cl}  the polygonal path $c_w$
is completely determined by the tree $\tilde{\theta_f} (f)$. 
It follows that $c_w (t)$ gives the same value for irrational  $w$ and $\t$
therefore  the polyhedron ${\cal N}(\psi_{f_i} (f_Y))$ 
is defined by the first formula of (\ref{New2}) by lemma \ref{curve}.

The case of the Newton polyhedron of $\psi_{f} (f_Y)$ is discussed 
analogously using that if $f_Y \in {\cal RC}_f$ and $f= b_1 \dots b_{s(f, f_Y)}$ is its bunch decomposition 
the result follows: the Minkowski sum  (\ref{New2}) corresponds to the decomposition 
$\psi_{f} (f_Y) = \prod_{i=1}^{s(f, f_Y)} \mbox{Res}_Y (T - f_1 \cdots f_s , b_i)$.
\hfill $\ {\diamond}$

\begin{Rem}
The Newton polyhedron  of $\psi_{f} (f_Y)$ is not necessarily polygonal since 
the set  $\{ v(P) \} $, for $P$ running through the non extremal vertices of
$\theta_f (f)$ is not totally ordered in general.  
\end{Rem}

\noindent
{\bf Proof of  Theorem \ref{main}.}
By Theorem   \ref{New} the Newton polyhedra 
$\psi_{f_i} (f_Y)$ for $i =1, \dots,s$, are polygonal 
and coincide with the those obtained assuming the hypothesis 
of $f_Y \in {\cal RC}_f$. Then assertion 1 follows  
by proposition \ref{Wall2}.

To prove assertion 2 
we consider the matrix ${\cal M} = (m_{i,j})$ whose columns are
the $t$-uples of
vectors defining the type of $f_Y$: 
$$ \left( \frac{\r(f_1, h_j)}{\deg{h_j}}, \dots,
\frac{\r(f_s, h_j)}{\deg{h_j}};   \deg {b_j} \right) \mbox{ for } {j=1, \dots, s(f,f_Y)}.$$
By definition, the columns of the matrix ${\cal M} $ correspond bijectively to the 
bunches of the decomposition of $f_Y$ induced by $f$. By theorem \ref{main} these bunches  correspond 
bijectively with the non extremal vertices of the tree $\theta_f (f)$, in such a way that 
if the column $j$ corresponds to the vertex $P_j$  then $m_{i, j} = \nu_i (P_j)$ for 
$i=1, \dots, s$.
We build the tree  $\theta_f (f)$ from the matrix ${\cal M}$ by  identifying the 
columns of  ${\cal M}$  with those non extremal vertices of $\theta_f (f_r)$
for $r=1, \dots, s$ separately:

We begin by analysing the row $r$.
If $a \in \{ m_{r, j} \} _{j=1} ^{s(f, f_Y)}$ the set of columns ${\cal K}_a ^r := \{ j / m_{r, j} = a \} $
of ${\cal M}$ is non empty and is clearly in bijection with the set 
${\cal P}_a ^r = \{ P \mbox{ vertex of } \theta_f (f) / \nu_r (P)  = a \}  $.
Since the set ${\cal P}_a ^r$ has a minimum for the valuation $v$, namely the vertex 
$Q $ of $\theta_f ( f_r)$ such that $\nu_r (P)  = a$, it follows by lemma \ref{monotone} that 
there is a unique column $l \in {\cal K}_a ^r$, corresponding to $Q$, 
such that $m_{t, l} \leq m_{t, k} $ for $t =1, \dots, s$ and $k \in {\cal K}_a ^r$. 
This procedure defines a partial order in the columns.
We recover the skeleton of the tree $\theta_f (f)$ by repeating this procedure for 
the rows $r=1, \dots,s $ of ${\cal M}$.
The vertex of bifurcation of $\theta_f (f_r)$ and  $\theta_f (f_k)$
is the greatest common column defined by the rows $r$ and $k$. 

To determine the chain $\g_f$ we use the row $s(f,f_Y) +1$ of ${\cal M}$ and 
the degrees of the irreducible factors of $f$. 
By theorem \ref{main} we know that the integer $m_{s(f,f_Y) +1, j}$ is the 
coefficient of the vertex $P_j$, corresponding to the column $j$, 
in the chain $-\partial \g_f$.  The coefficient of the extremal edge
containing the the vertex $P_{+\infty}^{(i)}$ in the chain $\g_f$
is equal to $\deg f_i$ for $i=1, \dots, s$. 

Since $\theta_f(f) $ is a tree, we recover recursively 
the coefficients appearing in the segments of the chain $\g_f$
from the chain 
$- \partial \g_f$
and the coefficients $\deg f_i$.
The chain $\g_f$ defines the vertices of $\theta_{f_i} (f_i) $ and the associated characteristic integers.
We recover from them, by using (\ref{19})
and the $\nu_i$ valuation,  the valuations $v(P)$ for those 
non extremal vertices  $P$ of $\theta_{f} (f_i) $.
\hfill $\ {\diamond}$

\section{A geometrical characterization of the bunch decomposition} \label{le}

We give a geometrical characterization of the bunch 
decomposition of $f_Y$
in terms 
of the {\it partial embedded resolution} $p : {\cal Z} \rightarrow \C^{d+1} $  of $f =0$ built in 
\cite{GP3} and  \cite{GP4}, by Gonz\'alez P\'erez. 
The morphism $p$  is a a composition of 
toric modifications which are 
canonically determined by the given quasi-ordinary projection, 
by using the tree ${\theta}_f (f)$.
An embedded resolution of $f=0$ is obtained by composing the modification $p$ with any
toroidal modification defining resolution of singularities of ${\cal Z}$, 
which always exists (see \cite{TE}).

The exceptional fiber $p^{-1} (0) $ of the modification $p$ is a curve, 
its irreducible components are 
complex projective lines.
The definition of the modification $p$ 
induces a  bijection $P \mapsto C(P)$
between the non extremal vertices of the tree $\theta_f (f)$ and the irreducible components of 
the exceptional fiber of $p$.  

Theorem  \ref{main} establishes a canonical bijection 
$P \rightarrow b_P$ between the non extremal vertices of the tree $\theta_f (f)$ 
and the bunches of the decomposition of $f_Y$ induced by $f$. 
An irreducible factor $h$ of $f_Y$ is a factor of 
of $b_P$ if and only if the following equality holds:
\begin{equation} \label{P}
\left( \frac{\r(f_1, h)}{\deg{h}}, \dots,
\frac{\r(f_s, h)}{\deg{h}} \right) = \left( \nu_1 (P), \dots, \nu_s (P)
 \right).
\end{equation}

Let  $h \in {\cal C}_f$ be irreducible, we say that 
$h$ is {\it associated} to a non extremal vertex $P$ of the tree  $\theta_f (f)$
if the equality (\ref{P}) holds, or equivalently
${\cal N}( \psi_{f_i} (h) ) =  \frac {\deg h}{\overline{\nu_i (P)}}$ for $i=1, \dots, s$;
this equivalence is deduced easily by arguing as in the proof of theorem \ref{New}. 
The following theorem implies that the $f$-bunch decomposition of $f_Y$ 
is compatible with the bijections above.

\begin{The} \label{lmw}
If $h \in {\cal C}_f$ is associated to a non extremal vertex $P$ of the tree  $\theta_f (f)$
then the strict transform of the hypersurface $h= 0$
only intersects the irreducible component  $C(P)$ of $p^{-1} (0) $.
The strict transform of $h= 0$  does not intersect 
the strict transform of $f =0$. 
\end{The}

We describe the procedure used to build the 
modification $p$ (for details see \cite{GP3}).
The modification $p$ is a composition $p  = \p_1 \dots \p_l$ of toric modifications.
The irreducible components of the exceptional fiber $p^{-1} (0)$ are complex projective lines
which can be ordered by $C' < C$ iff there exists $t>1$ such that  the image of
$C'$ by the modification $p': = \p_t \circ \dots \circ \p_l$ is a point of $p' (C)$ and 
$p'(C) $ is not reduced to a point. 
By definition the minimal components of $p^{-1} (0)$ with respect to this relation 
are the irreducible components of 
the exceptional fiber $\p_1 ^{-1} (0)$ of the first toric modification.

The morphism     $\p_1$ is the toric modification defined by the dual Newton diagram of 
$\Sigma (f) $ of the  quasi-ordinary polynomial $f \in \C \{ X \} [Y]$, when $Y$ is a {\it good coordinate}.
Such a good coordinate is built by a $ \C \{ X \}$-automorphism of the polynomial ring $ \C \{ X \} [Y]$
of the form $ Y \mapsto Y + r(X)$. These automorphisms are compatible with 
the sets ${\cal C}_f$ and ${\cal RC}_f$ since they preserve  resultants and discriminants 
of polynomials. 
We suppose from now on that $Y$ is a good coordinate for $f$, this means that 
the Newton polyhedron ${\cal N} (f) $ is polygonal
and it is completely determined from the Eggers-Wall tree $\theta_f (f)$.
In order  describe this polyhedron we define 
$${\cal A}_i^{(f)} :=    ( M \cap \{ k(f_i, f_j) \}_{j} ) \cup  \{ \l_1^{(i)} \}  \mbox{ for }     1 \leq i \leq s.$$  
By lemma \ref{rorder}, if the set  ${\cal A}_i^{(f)}$ is non empty it is totally ordered and we define then 
\begin{equation} \label{kappa}
\l_{\kappa (i)}^{(f)} := \left\{ \begin{array}{l}
 \min {\cal A}_i^{(f)} \mbox{ if }
{\cal A}_i^{(f)} \ne \emptyset
\\
{+ \infty} \mbox{ otherwise }
\end{array}\right\} \mbox{ for } i=1, \dots, s.
\end{equation}
The polyhedron  ${\cal N} (f_i)$ has only one compact edge (since is a Minkowski term of a polygonal polyhedron) and
since $f_i$ is irreducible it has only one compact edge of inclination equal to $ \l_{\kappa (i)}^{(f)}$ if 
 $\l_{\kappa (i)}^{(f)} \ne + \infty$ (the case $\l_{\kappa (i)} ^{(f)} = + \infty$ may happen
only for one index $i$ and in that case $f_i$ is a good coordinate for $f$, i.e., we will suppose that $f_i = Y$).

\begin{Lem} \label{Prepa} 
If the term $X^{\l}$  appears in the expansions 
of the roots of $f_j$ and if $\l_{\kappa (i)} ^{(f)} \nleq \l \notin M$ then we have that 
${\l} \geq k(f_i, f_j)$
and the  equality ${\l} = k(f_i, f_j) $ implies that $k(f_i, f_j)  = \l_1^{(j)}$.
The set
$\{ \l_{\kappa (1)}^{(f)}, \dots,  \l_{\kappa (s)} ^{(f)}\}$ is totally ordered by (\ref{order}) and  
its intersection with the reference lattice $M$ is defined by its maximal element or empty.
We have that 
$$k(f_i, f_j) <  \l_{\kappa (i)}^{(f)}  \Leftrightarrow \l_{\kappa (j)}^{(f)} < \l_{\kappa (i)}^{(f)}.$$
\end{Lem}
{\em Proof.}
For the assertions 1 and 2 see lemma 3.15 of \cite{GP3}.
Suppose that $k(f_i, f_j) <  \l_{\kappa (i)}^{(f)}$,
 if $k(f_i, f_j) \notin M$ it follows from assertion 1 that $k(f_i, f_j) = \l_1^{(j)}$.
In any case it follows from the definition that $\l_{\kappa (j)}^{(f)} < \l_{\kappa (i)}^{(f)}$.
Conversely, if $\l_{\kappa (j)}^{(f)} < \l_{\kappa (i)}^{(f)}$ then
 $\l_{\kappa (j)}^{(f)} \notin M$ by assertion 2 thus $\l_{\kappa (j)}^{(f)} = \l_1^{(j)}$ by definition, 
therefore $k(f_i, f_j) \leq \l_1^{(j)} < \l_{\kappa (i)}^{(f)}$ by assertion 1. 
\hfill $\ {\diamond}$

The exceptional fiber of the toric modification $\p_1$ is described by (\ref{exceptional}). 
We denote by $C ({P_{\l_{\kappa (i)}^{(f)}} ^{(i)}})$ the irreducible component of $\p_1^{-1} (0) $ which is  the closure of the orbit associated to 
the cone of $\Sigma(f)$ 
orthogonal to the compact edge of ${\cal N} (f)$
with inclination ${\l_{\kappa (i)}^{(f)}} \ne + \infty$.

The following lemma describes some properties of the strict transform 
of $f$ by $\p_1$.
\begin{Lem} \label{sep}
The strict transform of $f_i = 0 $ by $\p_1$ is a germ at the point of intersection 
$o_1^{(i)}$ with $\p_1^{-1} (0)$. This point belongs to only one irreducible component of 
 $\p_1^{-1} (0) $ which is equal to  $C( P^{(i)}_{ \l_{\kappa (i)}^{(f)}} )$
 if $ \l_{\kappa (i)}^{(f)} \ne + \infty$
or to $C( P^{(*)}_{\max \{ \l_{\kappa (j)}^{(f)} \}_{j=1}^s})$ otherwise. 
If $ \l_{\kappa (i)}^{(f)} \ne + \infty$ then we have that: 
$$C ({P_{ \l_{\kappa (i)}^{(f)}}^{(i)}} )  = C ({P_{ \l_{\kappa (j)}^{(f)}}^{(j)}} ) \Leftrightarrow 
k(f_i, f_j) \geq  \l_{\kappa (i)}^{(f)}
\mbox{ and in this case } 
o_1^{(i)} = o_1^{(j)} \Leftrightarrow \; k(f_i, f_j) > \l_{\kappa (i)}^{(f)} $$
\end{Lem}
{\em Proof.} 
See Proposition 3.32 of \cite{GP3}. The main assertion is also 
a direct consequence of the proof of theorem \ref{irred} and lemma \ref{Prepa}.
\hfill $\ {\diamond}$

\begin{Rem} \label{cle} $\,$ 
\begin{enumerate}
\item
The point $o_1^{(i)}$ is parametrized by the only root $c_{f_i}$ of the polynomial in one variable
$p(f_i, {\mathcal E}_i)$ defined from the symbolic 
restriction of $f$ to the compact edge ${\mathcal E}_i$ of
 ${\cal N} (f_i)$ by  (\ref{edge}).
\item 
If $ \l_{\kappa (i)}^{(f)} \ne + \infty$ and if 
$h \in \C \{ X \} [Y]$ is such that ${\cal N} (h) = \frac{ \deg h}{\overline{ \l_{\kappa (i)}^{(f)} }}$ and 
the polynomial $p(h, {\mathcal E})$ defined from the symbolic 
restriction of $h$ to the compact edge ${\mathcal E}$ of the polyhedron 
$\frac{ \deg h}{\overline{ \l_{\kappa (i)}^{(f)} }} $ by  (\ref{edge})
has only one root $c_h = c_{f_i}$, then the strict transform of $h=0$ by $\p_1$ is a 
germ at the point $o_1^{(i)}$ of $p_1^{-1} (0)$. 
\item
Suppose that $ \l_{\kappa (i)}^{(f)}=  \max \{  \l_{\kappa (l)}^{(f)} \} < + \infty$,
if $h \in \C \{ X \} [Y] $ is any polynomial such that  
${\cal N} (h) \subset \frac{ \deg h }{\overline{ {\l_{\kappa (i)}^{(f)}}}}$
and that the polynomial $p(h, {\mathcal E})$
has only one root  $c_{h} = 0 $, the strict transform of $h= 0$ by $\p_1$ 
only intersects  the component 
$C ({P_{ \l_{\kappa (i)}^{(f)}}^{(i)}} )$. 
\end{enumerate}
\end{Rem}
{\em Proof.} 
For the first and the second assertion of the following remark 
see the proof of theorem \ref{irred}.
The third assertion is obtained by the same argument used in the proof of lemma \ref{spe}. 
\hfill $\ {\diamond}$

The key inductive step in the embedded resolution procedure is that the 
strict transform of $f_i = 0 $ at the point $o_1^{(i)}$ is a {\it toric quasi-ordinary 
singularity}\footnote{This means that there is a finite projection onto a germ of affine toric 
variety which is unramified outside its torus (see \cite{GP3} and \cite{GP1})}, 
with a canonical ``quasi-ordinary'' projection.
These singularities are hypersurfaces of affine toric varieties, 
for instance the strict transform of $f = 0$ at the point $o_1^{(i)}$ are 
defined by the vanishing of a monic polynomial $f '$  in one variable with coefficients
in the ring $ \C \{ \D_d^\vee \cap (M + \l_{\kappa (i)}^{(f)} \Z)  \}$.
The polynomial $f'$ is {\it quasi-ordinary}: its discriminant is  
the product of a monomial by a unit of this ring. 
The definition of characteristic monomials and the Eggers-Wall tree 
introduced in the first section generalize to this setting by replacing the valuation 
$v$ by the lattice valuation $\tilde{v}$ and the 
reference lattice cone introduced in the first section 
(see  \cite{GP3}).
In particular the Eggers-Wall tree of $f'$ is determined from that of $f$
by the following proposition
(see Proposition 3.22 of \cite{GP3}):
\begin{Pro} \label{EW} If $ \l_{\kappa (i)}^{(f)} \ne + \infty$,  then we have that 
the Eggers-Wall tree $\theta_{f'} (f') $ 
associated to the strict transform of $f$ at the point $o_1^{(i)}$ is obtained from
$\theta_f (f) $ by removing  
the segment $ [ P_0^{(j)}, P^{(j)}_{\l_{\kappa (i)} ^{(f)}} [$
from 
the sub-tree of $\theta_f (f)$ given by $\bigcup \theta _f ( f_j )$, for 
those irreducible factors $f_j$ with order of coincidence $> \l_{\kappa (i)}^{(f)} $ with $f_i$.
The new lattice valuation is $\tilde{v}'(P) = \tilde{v} (P) - \l_{\kappa (i)}^{(f)}$.
The coefficients of the $1$-chain $\gamma_{f'} $ are obtained from those of $\gamma_{f} $
by division by the index of $\l_{\kappa (i)}^{(f)}$  over the old reference lattice $M$.
The new reference lattice cone is  $(\D_d^\vee, M + \l_{\kappa (i)}^{(f)} \Z)$.
\hfill $\ {\diamond}$
\end{Pro}

Proposition \ref{EW} allows us to extend the 
natural bijection between the components of the exceptional fiber of $\p_1$ and 
the subset $\{ P_{\l_{\kappa (1)}^{(f)}} ^{(1)}, \dots,    P_{\l_{\kappa (s)}^{(f)}} ^{(s)} \}$  
of vertices of $\theta_f (f)$, 
inductively between the 
the components of 
$p^{-1} (0)$ and the set of non extremal vertices of $ \theta_f (f)$.
The following lemma which translates some of the properties 
stated in terms of the notion of order of coincidence in terms of Newton polyhedra.
\begin{Lem} \label{cle2}
Let $h \in {\cal C}_f$ be an irreducible polynomial associated to the non extremal vertex $P$ 
of $\theta_f(f)$.  
Let  $C({P_{\l_{\kappa (i)}^{(f)}} ^{(i)}})$ be the unique component of $p^{-1} (0)$ 
which is  $\leq C(P)$. Then we have that ${\cal N} (h) \subset \frac{ \deg h }{\overline{ {\l_{\kappa (i)}^{(f)}}}}$
with equality if 
${\l_{\kappa (i)}^{(f)}} \ne \max \{ \l_{\kappa (j)}^{(f)} \}_{j =1} ^s$.
The polynomial $p(h, {\mathcal E}) \in \C[t]$ obtained from 
the symbolic 
restriction of $h$ to the compact edge ${\mathcal E}$ of 
the polyhedron $\frac{{ \deg h} }{\overline{ {\l_{\kappa (i)}^{(f)}}}}$ 
by (\ref{edge}) has only one complex root $c_h$ and we have that 
$c_h = c_{f_i} \Leftrightarrow v(P) >  \l_{\kappa (i)}^{(f)}$ 
(resp.
$c_h \ne c_{f_i} \Leftrightarrow v(P) =   \l_{\kappa (i)}^{(f)}$).
The case $c_h = 0 $ may  happen only if  $v(P) =  \max \{ \l_{\kappa (j)}^{(f)} \}_{j =1} ^s$.
\end{Lem}
{\em Proof.}
If $h \in {\cal RC}_f$,  we consider it as an irreducible factor of the quasi-ordinary
polynomial $fh$. 
The hypothesis means that the vertex $P_{k(f, h)}^h$ of $\theta_{fh} (f)$ 
belongs to $\theta_{f} (f)$, therefore  $\theta_{fh} (f) = \theta_{f} (f)$ and
$P = P_{k(f, h)}^h$ hence we have:
\begin{equation} \label{ieq2}
v(P) = k(f, h) = \max_{l=1, \dots, s} \{ k(f_l, h) \} = k(f_i, h) =  \l_{\kappa (i)}^{(f)}, \mbox{ and } 
\end{equation}
\begin{equation} \label{ieq}
\l_{\kappa (l)}^{(fh)} =  \l_{\kappa (l)}^{(f)} \mbox{ for } l =1, \dots, s.
\end{equation}
The exponents appearing  on the parametrizations of $f_i$ are $\geq \l_{\kappa (i)}^{(f)}$,
since $Y$ is a good coordinate for $f$. 
The equality $k(f_i, h) = \l_{\kappa (i)}^{(f)}$ implies that the same property holds 
for the parametrizations of $h$ and therefore   
${\cal N} (h) \subset \frac{ \deg h }{\overline{ {\l_{\kappa (i)}^{(f)}}}}$.
If these polyhedra are not equal the intersection of the compact face ${\mathcal E}$
with ${\cal N} (h)$ is reduced to the point $(0, \deg h)$ therefore 
$c_h = 0$.

These polyhedra are equal if and only if  $c_h \ne 0$. 
We show that this is always the case when 
$\l_{\kappa (i)}^{(f)} < \l_{\kappa (j)}^{(f)}$ for some index $j$: 
We deduce from lemma \ref{Prepa} 
and (\ref{ieq2})
that $ k(f_i, h) = \l_1^{(i)}$. 
We have that 
$k(h, f_j )  \leq  k(h, f_i ) <  \l_{\kappa (j)}^{(fh)} =  \l_{\kappa (j)}^{(f)}$
by (\ref{ieq2}) and (\ref{ieq}), therefore 
$k(h, f_j ) \notin M$ by definition of $\l_{\kappa (j)}^{(fh)}$.
By lemma \ref{Prepa} we deduce that $k(h, f_j ) = \l_1^{(h)} \leq  k(h, f_i ) = \l_1^{(i)}$
hence $\l_1^{(h)} = k(h, f_i )= \l_{\kappa (i)}^{(f)} $ by definition of order of coincidence. Therefore, we 
obtain that 
${\cal N} (h) = \frac{ \deg h }{\overline{ {\l_{\kappa (i)}^{(f)}}}}$.
The rest of the assertion in this case follows from lemma \ref{sep}.

We prove that this result holds also for $h \in {\cal C}_f$ by using  toric base changes  
of section \ref{tbc}. Let $\t \subset \D_d$ be any regular cone of dimension $d$ such that 
$h^{(\t)} \in {\cal RC}_{ f ^{(\t)}}$.  
The polynomial 
 $p(h, {\mathcal E})$ coincides with $p(h^{(\t)} , {\mathcal E})$ by lemmas \ref{soy} and \ref{tau2}.
The same argument provides the assertion about the inequality 
${\cal N} (h) \subset  \frac{ \deg h }{\overline{ {\l_{\kappa (i)}^{(f)}}}}$. 
\hfill $\ {\diamond}$

\noindent
{\bf Proof of Theorem \ref{lmw}.}
If $P = {P_{\l_{\kappa (i)}^{(f)}} ^{(i)}}$ for some $i$,  we deduce from lemma \ref{cle2} 
and remark \ref{cle} that 
the strict transform ${\cal H}'$  of $h = 0$ by $\p_1$ 
only meets the component $C(P)$ of $\p_1 ^{-1} (0)$ and it does not intersect 
the strict transforms of 
$f = 0$. This implies that a neighborhood of  ${\cal H}'$ wont be modified by the toric modifications 
$\p_2, \dots, \p_l$ and proves the theorem in this case.

If $C ({P_{\l_{\kappa (i)}^{(f)}} ^{(i)}}) <  C (P)$ we consider the 
polynomial $h'$ defining the strict transform of $h$ (which is a germ at the point $o_1^{(i)}$).
It follows from proposition  \ref{EW} and another application of lemmas \ref{soy} and \ref{tau2}
that $h'$ is associated to $P$ viewed on the tree $\theta_{f'} (f')$. 
If $C( P)$ is minimal between those components of $p^{-1} (0) $ corresponding 
to the non extremal vertices of $\theta_{f'} (f')$ we can apply the 
arguments of the previous case, otherwise we obtain the result by iterating the procedure. 
\hfill $\ {\diamond}$

\subsection{A theorem of L\^ e, Michel and Weber revisited} \label{slmw}

If $d= 1$, then $f \in \C \{ X \} [Y]$ defines the germ of a complex analytic plane curve 
in $(\C^2, 0)$.
The {\it minimal embedded resolution} of $(S, 0)$ is the modification 
$\Pi : {\cal X} \rightarrow \C^2 $  defined by the composition of the minimal sequence of 
points blow ups, such that the total transform 
of $f= 0$ is a normal crossing divisor. 
The {\it dual graph}  ${\cal G} (\Pi, 0) $ (resp. the {\it total dual graph}  ${\cal G} (\Pi, f) $) 
is the graph obtained from the exceptional divisor $\Pi ^{-1} (0) $ (resp. from 
the total transform  $\Pi ^{-1} ( \{ f =0 \}) $)
by associating a vertex to any irreducible component and joining with a segment 
those vertices whose associated components have non empty intersection. 
The {\it valency}  of 
a vertex $ P$ of a finite graph is the number of edges of the graph 
which contain the vertex $P$. 
We denote by $\# 1$ the component of  ${\cal G} (\Pi, f) $
which corresponds to the first blow up. 
We define for a vertex $P$ of ${\cal G} (\Pi, f) $ the integer 
\begin{equation} \label{omega}
\w(P) := 
 \left\{ \begin{array}{lcl}
  \mbox{valency of } P \mbox{ in   }    {\cal G} (\Pi, f)& \mbox{ if  } &  P \ne \# 1
\\
 1 + \mbox{ valency of } P \mbox{ in   } {\cal G} (\Pi, f) & \mbox{ if  } &  P = \# 1
\end{array}\right.
\end{equation}
A vertex $P$ of ${\cal G} (\Pi, f) $  with $\w(P) = 1$ and $P \ne \# 1$   (resp. $\w(P) \geq 3$)    
is called an  {\it extremal vertex} (resp. a {\it rupture vertex}). 
A {\it dead arc} of the graph ${\cal G} (\Pi, f) $ is a closed polygonal in ${\cal G} (\Pi, 0) $
joining a extremal vertex to any rupture vertex of ${\cal G} (\Pi, f) $, and which does not 
contain any other rupture vertex. 
The dual graph  ${\cal G} (\Pi, f) $ defines a natural stratification of
$\Pi^{-1} ( \{ f =0 \}) $, the $0$-dimensional strata are in bijection with the
segments of  ${\cal G} (\Pi, f) $, each segment corresponds to the intersection of the irreducible
components associated to the vertices. 
The $1$-dimensional strata are in bijection with the vertices, the stratum corresponding to 
a vertex of ${\cal G} (\Pi, f) $ is 
the set of points of the corresponding component which do not belong to any other component of 
 $\Pi^{-1} ( \{ f =0 \}) $. 
We can associate by this correspondence 
a subset of $\Pi^{-1} ( \{ f =0 \}) $ to any subgraph of  ${\cal G} (\Pi, f) $.

With these notations 
the main result of L\^e, Michel and Weber in  \cite{LMW} is: 
\begin{The}\label{lmw2} 
Denote by ${\cal Q} (\Pi, f) $ the set associated to 
the subgraph of  ${\cal G} (\Pi, f) $ defined by the rupture vertices and 
dead arcs. 
If $X=0$ is not contained in the tangent cone of 
the plane curve germ $f=0$ the intersection  of 
the exceptional divisor $\Pi^{-1} (0)$ with 
the strict transform of the polar curve $f_Y = 0$ 
is contained    ${\cal Q} (\Pi, f) $
and meets any connected component of  ${\cal Q} (\Pi, f) $.
\end{The} 

\noindent
{\em Proof}. If $X=0$ is not contained in the tangent cone of 
the plane curve germ $f=0$
the minimal embedded resolution  $\Pi$
is 
the composition of the partial embedded resolution $p: {\cal Z} \rightarrow \C^2$ 
used in the previous section, 
with a finite number 
of local toric modifications at the isolated  singular points of the normal variety ${\cal Z}$.
See \cite{GP3}, section 3.3.4 for details. 
The notions of {\it dual graph} ${\cal G} (p, 0) $  and {\it total dual graph} ${\cal G} (p, f)$ 
can be defined in an analogous way for $p$. 
In particular, we have that ${\cal G} (p, 0) $ is combinatorially isomorphic to the tree $\theta_f(f)$ minus 
its extremal segments and that there is a natural inclusion  
of the vertices of ${\cal G} (p, 0) $ 
in the vertices of  ${\cal G} ( \Pi, 0) $ whose image is the subset of rupture vertices of 
${\cal G} ( \Pi, f) $. The dual graph  ${\cal G} (p, f) $
is associated in an analogous manner 
to a natural stratification of $p^{-1} (\{  f = 0 \} )$, in such a way that 
we can associate by this correspondence 
a subset of $p^{-1} ( \{ f =0 \}) $ to any subgraph of  ${\cal G} (p, f) $.
We denote by ${\cal Q} (p, f)$ the subset of $ p^{-1} ( \{ f =0 \}) $ 
corresponding to the set of vertices of ${\cal G} (p, 0) $.

An irreducible  factor $h$ of $f_Y$ is associated to a non extremal  vertex $P$
of $\theta_f(f)$ by Theorem  \ref{main}.  The strict transform of $h$ by $p$ 
intersects $p^{-1} (0)$ at a smooth point $o_h$,  
which belongs to the component $C(P)$ associated to $P$ and 
does not intersect the strict transform of $f$, by theorem \ref{lmw}.
We deduce that the point $o_h$ belongs to ${\cal Q} (p, f) $ and that if 
 $\Pi = p' \circ p$ 
then 
$ {\cal Q} (\Pi, f) = (p')^{-1} (  {\cal Q} (p, f) )$ 
and the result follows.
\hfill $\ {\diamond}$

\begin{Rem}
If $C(P) $ is minimal for $p^{-1} (0)$ the component corresponding to  $C(P)$ in $\Pi ^{-1}(0)$ 
belongs to a dead arc if and only if $ v(P) =  \max \{ \l_{\kappa(j)}^{(f)} \} $, 
which is necessarily $< + \infty$. If $C(P)$ is not minimal, after some toric modifications, 
$C(P)$ is minimal for the strict transform of $f$ and 
an analogous result holds.
\end{Rem}
By definition there is no dead arc corresponding to $\min  \{ \l_{\kappa(j)}^{(f)} \}$
since the extremal point $\# 1$ is not considered for the definition of dead arcs.

\section{The case of Laurent quasi-ordinary polynomials}

The class of {\it Laurent quasi-ordinary polynomials} is introduced by Popescu-Pampu in  \cite{PP2}, see 
also \cite{PP1}, 
by analogy with the case of {\it meromorphic plane curves} studied by Abhyankar and Assi (see \cite{A-A}). 
He proves  a  decomposition theorem
for the derivative of any polynomial in the class 
such that the derivative itself is also Laurent quasi-ordinary. 
In this section we generalize this result to any Laurent quasi-ordinary polynomial
by translating
in an equivalent manner,  the properties of  the
bunch decomposition on the derivative of a quasi-ordinary polynomial 
from the  Laurent case to holomorphic case (characterized in Theorem \ref{main}).
This reduction is partially inspired  by an argument of 
Kuo and Parusinski 
comparing the plane curve meromorphic 
case to the 
holomorphic case (see \cite{KP}).

We denote by  $\C \langle X \rangle $  the ring of {\it Laurent power series} in $X =(X_1, \dots, X_d)$,
that is  
the ring of fractions $\C \{ X \} [ X_1 ^{-1}, \dots, X_d ^{-1}]$.
A  Laurent polynomial $F  \in \C \langle X \rangle [Y]$
admits a Newton polyhedron which is defined,  as usual in terms of the power series expansion 
$F = \sum c_{\a, i} X^\a Y^i$, as the convex hull of the set
$ \bigcup_{c_{\a, i} \ne 0} (\a, i) + (\R^{d}_{\geq 0}  \times \{ 0 \})$. 
A Laurent monic polynomial $F \in   \C  \langle X \rangle    [Y]$ is {\it quasi-ordinary} if the discriminant
$\D_Y (f)$ is of the form $\D_Y (f) = X^{\delta} \epsilon $ where 
$\delta \in \Z^d$ and 
$\epsilon $ is a unit 
in the ring of power series $\C \{ X \}$.
We can extend the definition of polynomials  comparable to  a quasi-ordinary polynomial to the 
Laurent case. 
In particular if  $F$ is a Laurent quasi-ordinary polynomial 
we define the type of $F_Y$ as  in the holomorphic case 
(see Definition
\ref{type}).
We relate Laurent monic polynomials with monic holomorphic polynomials by:
\begin{Lem}\label{L}
Let $F \in \C \langle X \rangle [Y] $ be a monic polynomial. Then there exists 
a vector
$q \in \Z^d$ such that the monic  polynomial $f$ defined by:
\begin{equation} \label{q}
f := X^{- \deg (F) \,  q}  \, F (X^q Y) 
\end{equation}
belongs to $\C \{ X \}  [Y] $.
In this case we have:
\begin{enumerate}
\item 
If $F$ is quasi-ordinary the same holds for $f$.
\item 
If $F = F_1 \cdots F_s$ is the factorization in irreducible monic polynomials 
the same holds for $f= f_1 \cdots f_s$  
where the polynomials $f_i$ are defined from $F_i$ by (\ref{q}).
\item
The polynomial $r$ defined 
from $R= F_Y$ by (\ref{q}) is equal to  $r = f_Y$. 
\end{enumerate}
\end{Lem}
{\em Proof.}
The polyhedron ${\cal N } (F) \subset \R^ d \times \R$ is contained in an affine cone 
$W$ of the form: $$W = [(0, \deg F) , (a, 0)] + (\R^ d \times \{ 0 \})$$
for some integral vector $a \in (\deg F) \Z^d $. 
If $q : = \frac{1}{\deg F} a$ then we obtain that the polynomial $f$ defined by (\ref{q})
belongs to $\C \{ X \}$. The idea is that 
the $1$-dimensional face $[(0, \deg F) , (a, 0)] $ of the cone $W$
corresponds to the segment $[(0, \deg f) , (0, 0)] $  in the Newton 
polyhedron of  $f$.
If $F = \sum_{i=0}^n a_i Y^i $ then $f = X^{ - \deg (F) q} (\sum_{i=0}^n a_i X^{i q} Y^i)$
 and if $R = \frac{\partial F}{\partial Y}$ then we obtain that  the polynomial 
$$\frac{\partial f}{\partial Y} = X^{ (- \deg (F)+1) q}  \left(  \sum_{i=1}^n i a_i X^{(i-1) q} Y^{i-1} \right) $$
is equal to $r$. 
It follows from quasi-homogeneity and homogeneity properties of the generic  discriminant
and resultant
(see \cite{GKZ}, pag. 398-399)  that 
$f$ is quasi-ordinary if $F$ is.
See also the proof of theorem 3 of \cite{GP1} for details.

It remains to prove assertion  2.
Recall that if a domain $A$ is integrally closed with fraction field $K$ 
then the factorization of a monic polynomial in $A[Y]$ 
as product of monic irreducible factors coincides over $A[Y] $ and $K[Y]$. 
It is easy to see that the rings $ \C \{ X \} \subset \C \langle X \rangle $ are integrally closed and 
have the same fraction field $L$. 
For any fixed $q \in \Z^d$ the multiplicative endomorphism of 
 $ \C \langle X \rangle [Y] \backslash \{ 0 \}$ 
defined by (\ref{q}) is an  automorphism
which preserves degrees and monic polynomials, and 
extends to a multiplicative autormorphism of $L [Y] \backslash \{ 0 \}$.
Therefore the factorization in monic irreducible factors of a monic polynomial $F \in  
\C \langle X \rangle [Y]$ corresponds by this mapping 
to the factorization  in monic irreducible factors of $f$ 
in $L [Y]$, if in addition $f $ belongs to $\C \{ X \}  [Y] $ then this factorization
holds over $\C\{ X \}  [Y] $ since $\C \{ X \}$ is integrally closed.
\hfill $\ {\diamond}$

Let $F $ be a Laurent quasi-ordinary polynomial.
The vector $q \in \Z^d$ will be that built in the proof of  lemma \ref{L}.
If 
$F_Y = H_1 \cdots H_r $ is the factorization in monic irreducible 
polynomials the same holds for $f_Y = h_1 \dots h_r $ by Lemma \ref{L} (the 
polynomials $h_i$ being defined from $H_i$ by (\ref{q})). 
We obtain, by using 
the quasi-homogeneity and homogeneity properties of the generic resultant (see \cite{GKZ} pag. 398-399), 
that for any $i, j$:  
$$\mbox{Res}_Y (f_j, h_i) 
=  X ^{ -\deg (H_i) \, \deg(F_j) \, q} \, \mbox{Res}_Y (F_j, H_i) $$
and since $\deg H_i = \deg h_i$  and $\deg F_i = \deg f_i$ for all $i$,   we obtain that:
\begin{equation} \label{T}
\frac{ \r(F_j, H_i)}{ \deg H_i} = \frac{ \r(f_j, h_i)}{ \deg h_i} + \deg(f_j) \, q.
\end{equation}
We have proved the following
result:
\begin{Pro} 
The bunches of the $F$-decomposition of $F_Y$ correspond to 
the bunches of the $f$-decomposition of $f_Y$ by the transformation (\ref{q}).
The type of $F_Y$ is obtained from the type of $f_Y$ by (\ref{T}). \hfill $\ {\diamond}$
\end{Pro}

\section{An example}

The polynomial  $f_{i, j} := (Y^2 - i X_1^3 X_2^2)^2 - j X_1^5 X_2^4 Y $ 
is  quasi-ordinary with characteristic exponents 
$\l_1 = (\frac{3}{2}, 1) $ and   $\l_2 = (\frac{7}{4}, \frac{3}{2})$, 
and integers $n_1 = n_2 =2$ for any $i, j \in \C^*$.
The equation is obtained by defining a deformation of the monomial variety 
associated to a quasi-ordinary hypersurface (see \cite{GP3}) in an analogous manner as
the deformation of the monomial curve associated to a plane branch studied in  \cite{Mono}.

The Eggers-Wall tree associated to the polynomial $f = f_{1, 1} f_{1, 2} f_{2,1} f_{2,2} $ is below.
We have that $v(P_1) = \l_1 $ and  $v(P_2) = v(P_3) = \l_2$.
The edges are labeled with the coefficients of the chain $\g_f$
thus we have that $-\partial \g_f = 4 P_1 + 6 P_2 +6 P_3$. 

\begin{figure}[htbp]
$$\epsfig{file=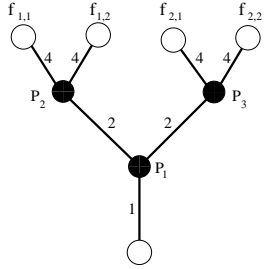, height= 35 mm}$$
\caption{The Eggers-Wall tree associated to $f$}
\end{figure}

We determine the type of $f_Y$ by using Proposition \ref{coincidence} and theorem 
\ref{main}:
\[
\begin{array}{cccccc}
f_{1,1} & \quad | & (6,4) & (\frac{13}{2}, 5) & (6,4) & | \\
f_{1,2} &\quad   |  & (6,4) & (\frac{13}{2}, 5) & (6,4) & |\\
f_{2,1} & \quad | & (6,4) &  (6,4) & (\frac{13}{2}, 5) & |\\
f_{2,2} & \quad|  & (6,4) &  (6,4) & (\frac{13}{2}, 5) & |\\
 & \quad | & 3 & 6 & 6 &|\\
& \quad   & - & - & - & \\
 & \quad & P_1 & P_2 & P_3 & \\
\end{array}
\]

\begin{figure}[htbp]
$$\epsfig{file=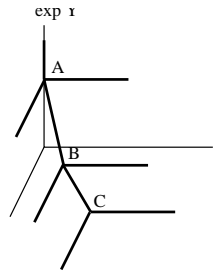, height= 35 mm}$$
\caption{A polygonal polyhedron \label{clp}}
\end{figure}

We compute the Newton polyhedra of the polynomials  $ \psi_{f_{1,1}} (f_Y)$ and $ \psi_{f} (f_Y)$ 
by using Theorem \ref{New}. 
We obtain that:
$
{\cal N} (\psi_{f_{i, j}} (f_Y) ) = \frac{3}{ \overline{(6,4)}} + \frac{6}{ \overline{(6,4)}} + 
\frac{6}{ \overline{( \frac{13}{2} , 5)}} = \frac{9}{ \overline{(6,4)}} + 
\frac{6}{ \overline{( \frac{13}{2} , 5)}}
$
for $i, j \in \{ 1, 2 \}$  and 
$
{\cal N} (\psi_{f} (f_Y) ) = \frac{3}{ \overline{(24 ,16)}} + \frac{6}{ \overline{(25, 18)}} + 
\frac{6}{ \overline{(25, 18)}} = \frac{3}{ \overline{(24 ,16)}} + 
\frac{12}{ \overline{(25, 18)}}$.
It follows that the polyhedra ${\cal N} (\psi_{f_{i, j}} (f_Y) )$ coincide in this example for $i, j \in \{ 1, 2 \}$.
In particular, the example shows that the only datum of these polyhedra does not allow us to distinguish
between the different irreducible factors of $f$. 
The polyhedron obtained for ${\cal N} (\psi_{f_{i, j}} (f_Y) )$ (resp. for  ${\cal N} (\psi_f (f_Y) )$)
is of the form of given in Figure \ref{clp} below where the vertices are 
$A = ( (0, 0), 15) $, $B = (( 54, 36), 6) $ and $C = ( ( 93, 66), 0)$;
(resp. $A = ( (0, 0), 15) $, $B = (( 72, 48), 12) $ and $C = ( ( 372, 264), 0)$).

If we are given the type of $f$ we recover the skeleton of the tree by noticing that
$\theta_f( f_{i, j}) $ has two non extremal vertices, corresponding to the different values 
appearing on the associated row. The corresponding minimal columns for $i=j =1$ (and for $i=1$, $j=2$ 
are $P_1 < P_2$.
Thus the column $P_2$ corresponds to the point of bifurcation of $\theta_f( f_{1, 1}) $
and $\theta_f( f_{1, 2}) $.
The columns corresponding to the non extremal vertices of $\theta_f( f_{i, j}) $ for $i=2, j=1, 2$
are $P_1 < P_3$, the bigger column corresponds with the point of bifurcation of $\theta_f( f_{2, 1}) $
and $\theta_f( f_{2, 2}) $ and the first is the point of bifurcation of $\theta_f( f_{1, 1}) $ and 
$\theta_f( f_{2, 1}) $.
The coefficient of the edge $\overline{P_1 P_2}$ on $\g_f$ is equal to
$\deg f_{1,1} + \deg f_{1, 2}  - c_{P_2} = 4+ 4- 6 = 2$ and we obtain the same value for 
$\overline{P_1 P_3}$. 
Then we recover the characteristic exponents $v(P_i)$ by using proposition 
\ref{coincidence}.

{ \small

\noindent
Evelia R. Garc\' \i a Barroso \\
Dpto. Matem\'atica Fundamental,
Fac. Matem\'aticas, Universidad de La Laguna\\
38271, La Laguna, Tenerife, Spain 

\noindent
{\tt ergarcia@ull.es}

\medskip

\noindent
Pedro Daniel Gonz\'alez P\'erez \\
Universit\'e de Paris 7, 
Institut de Math\'ematiques, 
Equipe G\'eom\'etrie et Dynamique\\
Case 7012;
2, Place Jussieu,
75005 Paris, France.

\noindent
{\tt gonzalez@math.jussieu.fr} }

\end{document}